\documentclass[a4paper,12pt]{article}
\usepackage[top=2.5cm,bottom=2.5cm,left=2.5cm,right=2.5cm]{geometry}
\usepackage{amsfonts}
\usepackage{latexsym,amsmath,amssymb,amsfonts,epsfig,psfrag,url,graphics,ifpdf,multicol}
\usepackage{color,colordvi,amscd,amsthm} %eepic
\usepackage{tikz}
\usepackage{verbatim}
\usepackage[numbers,sort&compress]{natbib}
\usepackage{indentfirst,graphics,epsfig}
\usepackage{graphicx}
\usepackage{graphics}
\usepackage{graphicx,psfrag}
\usepackage{graphics}
\usepackage{mathrsfs}
\usepackage{tabularx}
\usepackage{booktabs}
\usepackage{floatrow}
\usepackage{multirow}
\usepackage{graphicx}
\usepackage{caption}
\usepackage[ruled]{algorithm2e}
\newfloatcommand{capbtabbox}{table}[][\FBwidth]

\newtheorem{thm}{Theorem}
\newtheorem{lem}{Lemma}

\newtheorem{cl}{Claim}

\theoremstyle{plain}
\theoremstyle{plain}
\newtheorem{conj}{Conjecture}
 \theoremstyle{definition}

\theoremstyle{remark}
\def\pf{\noindent{\bf Proof.\ }}
\def\qed{{\hfill\rule{4pt}{7pt}}}

\allowdisplaybreaks

\numberwithin{subcase}{case}
\numberwithin{subca}{ca}
\numberwithin{subsubcase}{subcase}

\usepackage{indentfirst,subfig}
\begin{document}
	%\thispagestyle{empty}
	%\rule{0cm}{0.5mm}
	\captionsetup[figure]{labelfont={bf},name={Fig.},labelsep=period}

	\begin{center} {\large Upper bound for the number of maximal dissociation sets in trees}
	\end{center}
	\pagestyle{empty}
	
	\begin{center}
		{
			{\small Ziyuan Wang$^{1}$, Lei Zhang$^{2}$, Jianhua Tu$^{1,}$\footnote{Corresponding author.\\\indent \ \  E-mail: tujh81@163.com (J. Tu)}, Liming Xiong$^2$}\\[2mm]
			
			{\small $^{1}$School of Mathematics and Statistics, Beijing Technology and Business University, \\
				\hspace*{1pt} Beijing, P.R. China 100048}\\
			{\small $^2$ School of Mathematics and Statistics, Beijing Institute of Technology, \\
				\hspace*{1pt} Beijing, P.R. China 100081} \\[2mm]}
		
	\end{center}
	
	\begin{center}
		\begin{abstract}
			Let $G$ be a simple graph. A dissociation set of $G$ is defined as a set of vertices that induces a subgraph in which every vertex has a degree of at most 1. A dissociation set is maximal if it is not contained as a proper subset in any other dissociation set. We introduce the notation $\Phi(G)$ to represent the number of maximal dissociation sets in $G$. This study focuses on trees, specifically showing that for any tree $T$ of order $n\geq4$, the following inequality holds:			
	    	\[\Phi(T)\leq 3^{\frac{n-1}{3}}+\frac{n-1}{3}.\]			
			We also identify the extremal tree that attains this upper bound. Additionally, to establish the upper bound on the number of maximal dissociation sets in trees of order $n$, we also determine the second largest number of maximal dissociation sets in forests of order $n$.
			
				\vskip 3mm
			\noindent\textbf{Keywords:} Maximal dissociation sets; Enumeration; Trees; Forests\\
			\noindent\textbf{Mathematics Subject Classification:} 05C30; 05C69; 05C05 
			
		\end{abstract}
		%\end{minipage}
	\end{center}

	\baselineskip=0.24in
	%---------------------------------------
	\section{Introduction}\label{sec1}
Throughout this paper, we consider only simple, undirected, and labelled graphs. The order of a graph is defined as the number of vertices it contains.
Let $P_n$ represent the path and $K_{1,n-1}$ represent the star, both on $n$ vertices. Given two graphs $G$ and $H$, we define their disjoint union as $G\cup H$. Furthermore, we denote the disjoint union of $k$ copies of $G$ as $kG$.

%Unless otherwise stated, we follow the traditional notation and terminology (see also \cite{Bondy2008}).

An independent set in a graph is defined as a set of vertices in which no two vertices are adjacent to each other. Around 1960, Erd\H{o}s and Moser posed the intriguing questions:

{\it What is the maximum possible number of maximal independent sets in general graphs of order $n$? And what are the graphs that achieve this maximum number? }

In 1965, Moon and Moser \cite{Moon1965} provided answers to these queries, sparking a renewed interest in exploring the problem of determining the maximum number of maximal or maximum independent sets in special families of graphs. This pursuit also involves characterizing the graphs that achieve this maximum number. This type of problem is also referred to as the extremal enumeration problem.

Regarding the extremal enumeration problem for maximal or maximum independent sets, various families of graphs have been examined, including trees, forests, connected graphs, bipartite graphs, unicyclic graphs, graphs with at most $r$ cycles, and others. For a comprehensive overview of these findings, we refer the reader to \cite{Furedi1987, Griggs1988, Koh2008, Liu1994, Sagan1988, Sagan2006, Taleskii2022, Wilf1986, Ying2006, Zito1991}. Mohr and Rautenbach \cite{Mohr2018,Mohr2021} delved into connected graphs and trees of given order and independence number, while Taletskii and Malyshev \cite{Taleskii2022} focused on trees with a given number of leaves.

For the extremal enumeration problem, researchers have delved into various other graph substructures, including independent sets \cite{Davies2017, Sah2019}, matchings \cite{Davies2017}, maximal or maximum matchings \cite{Gorska2007, Heuberger2011}, maximal induced matchings \cite{Basavaraju2016}, minimal or minimum dominating sets \cite{Alvarado2019, Couturier2013}, and $k$-dominating independent sets \cite{Gerbner2019}, among others.

Let $G$ be a graph. A dissociation set of $G$ refers to a set of vertices that induces a subgraph composed of isolated vertices and paths $P_2$. A dissociation set is maximal if it is not a proper subset of any other dissociation set, and maximum if it contains the largest possible number of vertices. The dissociation number of $G$ represents the cardinality of a maximum dissociation set in $G$.

The concept of dissociation set was introduced by Yannakakis \cite{Yannakakis1981} in the early 1980s. In the past four decades, researchers have examined the dissociation set from various perspectives. For a comprehensive overview of related findings, we refer the reader to \cite{Bock2023-1, Bock2023-2, Cheng2023, Orlovich2011, Sun2023, Tu2021, Yannakakis1981}. Notably, recent years have seen a surge of interest in the extremal enumeration problem regarding maximal or maximum dissociation sets. Tu, Zhang, and Shi \cite{Tu2021} established that for any tree $T$ of order $n$, the number of {\bf maximum dissociation sets} in $T$ is upper-bounded by

\[\left\{
\begin{array}{ll}
	3^{\frac{n}{3}-1}+\frac{n}{3}+1, & \hbox{if $n\equiv0\pmod{3}$,} \\
	3^{\frac{n-1}{3}-1}+1, & \hbox{if $n\equiv1\pmod{3}$,} \\
	3^{\frac{n-2}{3}-1}, & \hbox{if $n\equiv2\pmod{3}$,}
\end{array}
\right.
\]
and characterized the extremal trees that achieve these upper bounds. Tu, Li, and Du \cite{Tu2022} provided tight upper bounds for the number of maximal and maximum dissociation sets in general graphs and triangle-free graphs of order $n$. More recently, Sun and Li \cite{Sun2023} determined the largest number of maximum dissociation sets in forests of given order and dissociation number.

We write $\Phi(G)$ to denote the number of maximal dissociation sets in a graph $G$.
Let \begin{align*}
	\begin{split}
		f_1(n):=\left \{
		\begin{array}{ll}
			3^{\frac{n}{3}},              & \hbox{if $n\equiv0\pmod{3}$ and $n\geq 3$,} \\
			4\cdot 3^{\frac{n-4}{3}},     & \hbox{if $n\equiv1\pmod{3}$ and $n\geq 4$, }\\
			5,                            & \hbox{if $n=5$,} \\
			4^2\cdot 3^{\frac{n-8}{3}},   & \hbox{if $n\equiv2\pmod{3}$ and $n\geq 8$.}
		\end{array}
		\right.
	\end{split}
\end{align*}
It is worth noting that for any integer $n\geq3$, $f_1(n)\leq 3^{\frac{n}{3}}$. Recently, Cheng and Wu \cite{Cheng2023} determined the largest number of maximal dissociation sets in forests of order $n$. They gave the following results. 
	
\begin{thm}\label{thm1}\cite{Cheng2023}\label{thm1}
		For any forest $F$ of order $n\geq 3$, \[\Phi(F)\leq f_1(n),\] with equality if and only if
			\begin{align*}
			\begin{split}
				F\cong \left \{
				\begin{array}{ll}
					\frac{n}{3}P_3,              & \hbox{if $n\equiv0\pmod{3}$,} \\
					\frac{n-4}{3}P_3\cup K_{1,3},  & \hbox{if $n\equiv1\pmod{3}$,} \\
					K_{1,4},                            & \hbox{if $n=5$,} \\
					\frac{n-8}{3}P_3\cup 2K_{1,3},   & \hbox{if $n\equiv2\pmod{3}$\ and $n\geq 8$.}
				\end{array}
				\right.
			\end{split}
		\end{align*}
	\end{thm}

Based on the results mentioned above, the following interesting and meaningful questions are raised.

{\it What is the largest number of {\bf maximal dissociation sets} in trees of order $n$? What are the extremal trees attaining the largest number?}

 In this paper, we focus on the above questions and give an upper bound on the number of maximal dissociation sets in trees of order $n$. We also characterize the extremal tree attaining the upper bound, thus showing that the upper bound is tight.

When $n\geq4$ and $n\equiv 1\pmod3$, let $T^*_n$ be the tree of order $n$ obtained from $\frac{n-1}{3}P_3$ and an isolated vertex by adding $\frac{n-1}{3}$ edges between the isolated vertex and the $\frac{n-1}{3}$ non-leaf vertices of $\frac{n-1}{3}P_3$. See Figure \ref{fig1}.

		\begin{figure}[H]
	\begin{tikzpicture}
		[line width = 1pt, scale=0.8,
		empty/.style = {circle, draw, fill = white, inner sep=0mm, minimum size=2mm}, full/.style = {circle, draw, fill = black, inner sep=0mm, minimum size=2mm}, full1/.style = {circle, draw, fill = black, inner sep=0mm, minimum size=0.5mm}]
		\node [empty] (a) at (0,2.4) {};
		\node [empty] (b) at (-2.4,1) {};
		\node [empty] (c) at (-2.7,0) {};
		\node [empty] (d) at (-2.1,0) {};
		\node [empty] (e) at (-1,1) {};
		\node [empty] (f) at (-1.3,0) {};
		\node [empty] (g) at (-0.7,0) {};
		\node [full1] at (0.8,0.7) {};
		\node [full1] at (0.4,0.7) {};
		\node [full1] at (1.2,0.7) {};
		\node [empty] (h) at (2.7,1) {};
		\node [empty] (i) at (2.4,0) {};
		\node [empty] (j) at (3.0,0) {};
		\draw (a) -- (b);
		\draw (a) -- (e);
		\draw (a) -- (h);
		\draw (b) -- (c);
		\draw (b) -- (d);
		\draw (e) -- (f);
		\draw (e) -- (g);
		\draw (h) -- (i);
		\draw (h) -- (j);
	\end{tikzpicture}
	\caption{$T^*_n$}
	\label{fig1}
\end{figure}

Let $T^*_8$ be the tree of order $8$ obtained from $2K_{1,3}$ by adding an edge between two leaves of $2K_{1,3}$. See Figure \ref{fig2}. It will be shown in Lemma \ref{lem4} that $T^*_8$ has the most number of maximal dissociation sets among all trees of order $8$.

\begin{figure}[H]
	\begin{tikzpicture}
		[line width = 1pt, scale=0.8,
		empty/.style = {circle, draw, fill = white, inner sep=0mm, minimum size=2mm}, full/.style = {circle, draw, fill = black, inner sep=0mm, minimum size=2mm}, full1/.style = {circle, draw, fill = black, inner sep=0mm, minimum size=0.5mm}]
		
		\node [empty] (h) at (5.5,0.5) {};
		\node [empty] (i) at (6.5,0.5) {};
		\node [empty] (j) at (6,0) {};
		\node [empty] (k) at (6,-0.5) {};
		\node [empty] (l) at (6,-1) {};
		\node [empty] (m) at (6,-1.5) {};
		\node [empty] (n) at (5.5,-2) {};
		\node [empty] (o) at (6.5,-2) {};
		\draw (h) -- (j);
		\draw (i) -- (j);
		\draw (j) -- (k);
		\draw (l) -- (k);
		\draw (l) -- (m);
		\draw (n) -- (m);
		\draw (o) -- (m);
		%\path (6,-2.5) node[below]{$T_{8max}$};
	\end{tikzpicture}
	\caption{$T^*_{8}$}
	\label{fig2}
\end{figure}

Let \[t(n):= 3^{\frac{n-1}{3}}+\frac{n-1}{3},\] and 
\begin{align*}
	\begin{split}
		f_2(n):=\left \{
		\begin{array}{ll}
			3 &\hbox{if $n=4$,}\\
			4 &\hbox{if $n=5$,}\\
			6 &\hbox{if $n=6$,}\\
			11\cdot 3^{\frac{n-7}{3}},     & \hbox{if $n\equiv1\pmod{3}$ and $n\geq 7$,} \\
			15\cdot 3^{\frac{n-8}{3}},     & \hbox{if $n\equiv2\pmod{3}$ and $n\geq 8$,} \\
			20,                            & \hbox{if $n=9$,} \\
			4^3\cdot 3^{\frac{n-12}{3}},   & \hbox{if $n\equiv0\pmod{3}$ and $n\geq 12$.}
		\end{array}
		\right.
	\end{split}
\end{align*}

The main contributions of our work are as follows.	
	
\begin{thm}\label{thm2}
	
(1) For any tree $T$ of order $n\geq4$, \[\Phi(T)\leq t(n),\]
with equality if and only if $n\equiv 1\pmod{3}$ and $T\cong T^*_n$.

(2) For any forest $F$ of order $4\leq n\leq 6$, if $\Phi(F)<f_1(n)$, then \[\Phi(F)\leq f_2(n).\] 
For any forest $F$ of order $n\geq7$, if $\Phi(F)<f_1(n)$, then \[\Phi(F)\leq f_2(n),\]
with equality if and only if
		\begin{align*}
			\begin{split}
				F\cong \left \{
				\begin{array}{ll}
					\frac{n-7}{3}P_3\cup T^*_{7}, & \hbox{if $n\equiv1\pmod{3}$,} \\
					\frac{n-8}{3}P_3\cup T^*_{8}\ or\ \frac{n-5}{3}P_3\cup K_{1,4},     & \hbox{if $n\equiv2\pmod{3}$,} \\
					K_{1,3}\cup K_{1,4},                            & \hbox{if\ $n=9$,} \\
					\frac{n-12}{3}P_3\cup 3K_{1,3},   & \hbox{if $n\equiv0\pmod{3}$ and $n\geq 12$.}
				\end{array}
				\right.
			\end{split}
		\end{align*}
\end{thm}

The rest of this paper is organized as follows. In Section \ref{sec2} we present four preliminary lemmas that play important role in the proof of Theorem \ref{thm2} which is given in Section \ref{sec3}. Some of the details of the proof of Theorem \ref{thm2} are included in Appendix. Finally, we conclude this paper in Section \ref{sec4}.
	
\section{Preliminary lemmas}\label{sec2}

Let $G$ be a graph, and define its {\it diameter} as the length of a longest path in it. For a vertex $v$ in $G$, we denote the set of its neighbors as $N_G(v)$ or simply $N(v)$. The {\it degree} of $v$ is denoted by $d(v)$ and represents the number of edges incident to $v$. A vertex with $d(v)=1$ is referred to as a {\it leaf} of $G$. A vertex adjacent to a leaf is termed a {\it support vertex}. The set of all maximal dissociation sets in $G$ is denoted by $MD(G)$. For any subset $A$ of the vertex set $V(G)$, we define
\[\Phi_A(G)=|\{S\in MD(G):A\subset S\}|,\ \Phi_{\overline{A}}(G)=|\{S\in MD(G):A\cap S=\emptyset\}|.\]
In particular, when $A$ consists of a single vertex $v$, we simplify the notation as
\[\Phi_v(G)=|\{S\in MD(G):v\in S\}|,\ \Phi_{\overline{v}}(G)=|\{S\in MD(G):v\notin S\}|.\]

\begin{lem}\label{lem1}
Let $T$ be a tree. If there exists a support vertex $v$ of degree $2$ in $T$ that is adjacent to a leaf $u$ and a non-leaf vertex $x$, then 
\[\Phi(T)\leq \Phi(T'),\]
where $T'$ is the tree obtained from $T$ by deleting the edge $uv$ and adding a new edge $ux$, i.e., $T'=T-uv+ux$. See Figure \ref{lem1-fig1}. Moreover, if  $N_T(x)=\{v,y_1,\cdots,y_{\ell}\}$ and there exists a component in $T-\{u,v,x,y_1,\cdots,y_{\ell}\}$ of order at least 3, then 
\[\Phi(T)<\Phi(T').\]
	
\begin{figure}[H]
	\begin{tikzpicture}
				[line width = 1pt, scale=0.8,
				empty/.style = {circle, draw, fill = white, inner sep=0mm, minimum size=2mm}, full/.style = {circle, draw, fill = black, inner sep=0mm, minimum size=2mm}, full1/.style = {circle, draw, fill = black, inner sep=0mm, minimum size=0.5mm}]
				\node [empty, label = right:$u$] (a) at (-3,2) {};
				\node [empty, label = right:$v$] (b) at (-3,1.5) {};
				\node [empty, label = right:$x$] (c) at (-3,1) {};
				\node [empty, label = right:$y_{\ell}$] (d) at (-2,0) {};
				\node [empty, label = right:$y_2$] (e) at (-3.5,0) {};
				\node [full1] at (-2.6,0) {};
				\node [full1] at (-2.45,0) {};
				\node [full1] at (-2.3,0) {};
				\node [empty, label = left:$y_1$] (f) at (-4,0) {};
				\draw (a) -- (b);
				\draw (b) -- (c);
				\draw (c) -- (d);
				\draw (c) -- (e);
				\draw (c) -- (f);
				\draw (-4.1,-1) -- (f);
				\draw (-3.9,-1) -- (f);
				\draw (-3.5,-1) -- (e);
				\draw (-2.3,-1) -- (d);
				\draw (-2,-1) -- (d);
				\draw (-1.7,-1) -- (d);
				\draw[black][dashed] (-5,1) rectangle (-1,-1.2);
				
				\draw[->](-0.75,0) -- (-0.25,0);
				\path (-3.25,-1.7) node(text1)[right]{$T$};
				
				\node [empty, label = above:$u$] (a) at (1.5,1.5) {};
				\node [empty, label = above:$v$] (b) at (2.5,1.5) {};
				\node [empty, label = above:$x$] (c) at (2,1) {};
				\node [empty, label = left:$y_1$] (d) at (1,0) {};
				\node [empty, label = right:$y_2$] (e) at (1.5,0) {};
				\node [full1] at (2.4,0) {};
				\node [full1] at (2.55,0) {};
				\node [full1] at (2.7,0) {};
				\node [empty, label = right:$y_{\ell}$] (f) at (3,0) {};
				\draw (a) -- (c);
				\draw (b) -- (c);
				\draw (c) -- (d);
				\draw (c) -- (e);
				\draw (c) -- (f);
				\draw (1.1,-1) -- (d);
				\draw (0.9,-1) -- (d);
				\draw (1.5,-1) -- (e);
				\draw (3.3,-1) -- (f);
				\draw (3,-1) -- (f);
				\draw (2.7,-1) -- (f);
				\draw[black][dashed] (0,1) rectangle (4,-1.2);
				\path (1.75,-1.7) node(text1)[right]{$T'$};
	\end{tikzpicture}
		\caption{$T'=T-uv+ux$}
		\label{lem1-fig1}
	\end{figure}
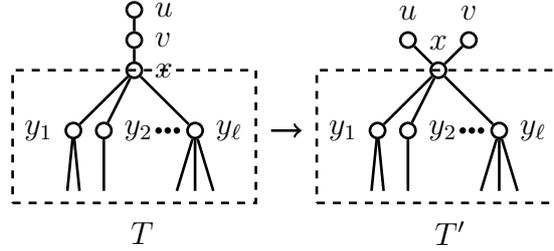
\end{lem}

\pf We have
	\[\Phi(T)=\Phi_{\{u,v\}}(T)+\Phi_{\{v,x\}}(T)+\sum\limits_{i=1}^{\ell}\Phi_{\{x,y_i\}}(T)+\Phi_{\overline{N(x)}}(T)\]
	and
	\begin{align*}
		\Phi(T')&=\Phi_{\{u,v\}}(T')+\Phi_{\{v,x\}}(T')+\sum\limits_{i=1}^{\ell}\Phi_{\{x,y_i\}}(T')+\Phi_{\{u,x\}}(T')\\
		&=\Phi_{\{u,v\}}(T)+\Phi_{\{v,x\}}(T)+\sum\limits_{i=1}^{\ell}\Phi_{\{x,y_i\}}(T)+\Phi_{\{u,x\}}(T').
	\end{align*}
If $S$ is a maximal dissociation set of $T$ such that $S\cap N(x)=\emptyset$, then $\{u,x\}\subseteq S$ and for each vertex $y_i\in\{y_1,\cdots,y_{\ell}\}$, there must be a vertex, say $z_i$, in $T-\{u,v,x,y_1,...,y_{\ell}\}$ such that $y_iz_i\in E(T)$ and $z_i\in S$. Thus, $\Phi_{\overline{N(x)}}(T)\leq \Phi_{\{u,x\}}(T')$, which implies that $\Phi(T)\leq\Phi(T')$. Furthermore, if there exists a component in $T-\{u,v,x,y_1,...,y_{\ell}\}$ of order at least $3$, then $\Phi_{\overline{N(x)}}(T)< \Phi_{\{u,x\}}(T')$, it follows that $\Phi(T)<\Phi(T')$.\qed

\begin{lem}\label{lem2}
	Let $T$ be a tree and $k\geq2$ be an integer. Let $x$ be a support vertex in $T$ of degree $k+1$ which is adjacent to $k$ leaves $u_1,\cdots,u_k$ and a non-leaf vertex $t$. If the vertex $t$ is adjacent to a leaf $y$, then 
	\[\Phi(T)\leq\Phi(T'),\]
	where $T'=T-xt+xy$. See Figure \ref{lem2-fig1}. Moreover, if $|V(T)\setminus\{u_1,\cdots,u_k,x,y\}|\geq3$, then \[\Phi(T)<\Phi(T').\]
	\begin{figure}[H]
		\begin{tikzpicture}
				[line width = 1pt, scale=0.9,
				empty/.style = {circle, draw, fill = white, inner sep=0mm, minimum size=2mm}, full/.style = {circle, draw, fill = black, inner sep=0mm, minimum size=2mm}, full1/.style = {circle, draw, fill = black, inner sep=0mm, minimum size=0.5mm}]
				\node [empty, label = above:$u_1$] (a) at (-4,3) {};
				\node [empty, label = above:$u_2$] (b) at (-3.5,3) {};
				\node [full1] at (-3,3) {};
				\node [full1] at (-2.75,3) {};
				\node [full1] at (-2.5,3) {};
				\node [empty, label = above:$u_k$] (r) at (-2,3) {};
				\node [empty, label = right:$x$] (c) at (-3,2) {};
				\node [empty, label = right:$y$] (d) at (-2.5,1.5) {};
				\node [empty, label = right:$t$] (e) at (-3,1) {};
				\node [empty, label = right:$y_{\ell}$] (f) at (-2,0) {};
				\node [empty, label = below:$y_2$] (g) at (-3.5,0) {};
				\node [empty, label = left:$y_1$] (h) at (-4,0) {};
				\node [full1] at (-3.0,0) {};
				\node [full1] at (-2.75,0) {};
				\node [full1] at (-2.5,0) {};
				\draw (a) -- (c);
				\draw (b) -- (c);
				\draw (c) -- (e);
				\draw (d) -- (e);
				\draw (e) -- (f);
				\draw (e) -- (g);
				\draw (e) -- (h);
				\draw (r) -- (c);
				\draw (-3.9,-1) -- (h);
				\draw (-4.1,-1) -- (h);
				\draw (-3.5,-1) -- (g);
				\draw (-1.8,-1) -- (f);
				\draw (-2.2,-1) -- (f);
				\draw (-2.0,-1) -- (f);
				\draw[black][dashed] (-5,1) rectangle (-1,-1.2);
				
				\draw[->](-0.75,0) -- (-0.25,0);
				\path (-3.25,-1.7) node(text1)[right]{$T$};
				
				\node [empty, label = above:$u_1$] (i) at (1,3) {};
				\node [empty, label = above:$u_2$] (j) at (1.5,3) {};
				\node [full1] at (2,3) {};
				\node [full1] at (2.25,3) {};
				\node [full1] at (2.5,3) {};
				\node [empty, label = above:$u_k$] (q) at (3,3) {};
				\node [empty, label = right:$x$] (k) at (2,2) {};
				\node [empty, label = right:$y$] (l) at (2,1.5) {};
				\node [empty, label = right:$t$] (m) at (2,1) {};
				\node [empty, label = left:$y_1$] (n) at (1,0) {};
				\node [empty, label = below:$y_2$] (o) at (1.5,0) {};
				\node [full1] at (2.0,0) {};
				\node [full1] at (2.25,0) {};
				\node [full1] at (2.5,0) {};
				\node [empty, label = right:$y_{\ell}$] (p) at (3,0) {};
				\draw (i) -- (k);
				\draw (j) -- (k);
				\draw (k) -- (l);
				\draw (m) -- (l);
				\draw (m) -- (n);
				\draw (m) -- (o);
				\draw (m) -- (p);
				\draw (q) -- (k);
				\draw (n) -- (0.9,-1);
				\draw (n) -- (1.1,-1);
				\draw (o) -- (1.5,-1);
				\draw (p) -- (2.8,-1);
				\draw (p) -- (3.0,-1);
				\draw (p) -- (3.2,-1);
				\draw[black][dashed] (0,1) rectangle (4,-1.2);
				\path (1.75,-1.7) node(text1)[right]{$T'$};
		\end{tikzpicture}
	   \caption{$T'=T-xt+xy$}
		\label{lem2-fig1}
	\end{figure}
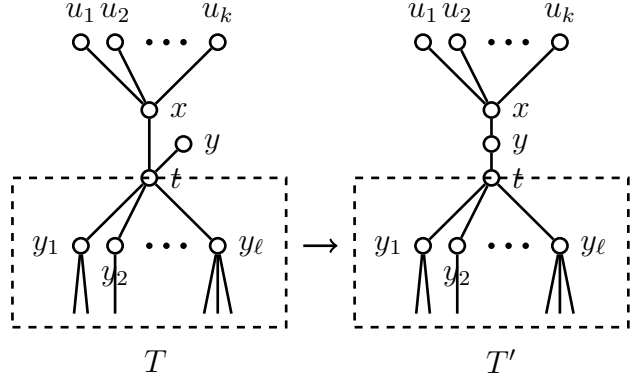
\end{lem}

\pf It is easy to see that $\Phi_{\{x,t\}}(T)\leq\Phi_{\{x,y\}}(T')$ and for any $i\in\{1,\cdots,k\}$, $\Phi_{\{u_i,x\}}(T)\leq\Phi_{\{u_i,x\}}(T')$. Thus, we have
		\[\Phi(T)=\sum\limits_{i=1}^k\Phi_{\{u_i,x\}}(T)+\Phi_{\{x,t\}}(T)+\Phi_{\{u_1,\cdots,u_k\}}(T),\]
		and
		\begin{align*}
		\Phi(T')&=\sum\limits_{i=1}^k\Phi_{\{u_i,x\}}(T')+\Phi_{\{x,y\}}(T')+\Phi_{\{u_1,\cdots,u_k\}}(T')\\
		&\geq\sum\limits_{i=1}^k\Phi_{\{u_i,x\}}(T)+\Phi_{\{x,t\}}(T)+\Phi_{\{u_1,\cdots,u_k\}}(T)\\
		&=\Phi(T).
		\end{align*}
	Moreover, if $|V(T)\setminus\{u_1,\cdots,u_k,x,y\}|\geq3$, then \[\sum\limits_{i=1}^k\Phi_{\{u_i,x\}}(T)+\Phi_{\{x,t\}}(T)<\sum\limits_{i=1}^k\Phi_{\{u_i,x\}}(T')+\Phi_{\{x,y\}}(T'),\]
	which implies that $\Phi(T)<\Phi(T')$. \qed

\begin{lem}\label{lem3}
Let $T$ be a tree of order $n\geq5$ and $v$ be a support vertex in $T$ that is adjacent to at least three leaves. Let $T'$ be a tree of order $n-1$ obtained from $T$ by deleting a leaf of $T$ that is adjacent to the vertex $v$. If $\Phi(T')\leq t(n-1)$, then $\Phi(T)<t(n).$
\end{lem}

\begin{figure}[H]
	\begin{tikzpicture}
		[line width = 1pt, scale=0.8,
		empty/.style = {circle, draw, fill = white, inner sep=0mm, minimum size=2mm}, full/.style = {circle, draw, fill = black, inner sep=0mm, minimum size=2mm}, full1/.style = {circle, draw, fill = black, inner sep=0mm, minimum size=0.5mm}]
		\node [empty, label = above:$u_1$] (a) at (1.0,2) {};
		\node [empty, label = above:$u_2$] (o) at (1.5,2) {};
		\node [full1] at (2.0,2) {};
		\node [full1] at (2.25,2) {};
		\node [full1] at (2.5,2) {};
		\node [empty, label = above:$u_k$] (b) at (3,2) {};
		\node [empty, label = above:$v$] (c) at (2,1) {};
		\node [empty, label = left:$y_1$] (d) at (1,0) {};
		\node [empty, label = right:$y_2$] (e) at (1.5,0) {};
		\node [full1] at (2.4,0) {};
		\node [full1] at (2.55,0) {};
		\node [full1] at (2.7,0) {};
		\node [empty, label = right:$y_{\ell}$] (f) at (3,0) {};
		\draw (a) -- (c);
		\draw (b) -- (c);
		\draw (c) -- (d);
		\draw (c) -- (e);
		\draw (c) -- (f);
		\draw (c) -- (o);
		\draw (d) -- (0.9,-1);
		\draw (d) -- (1.1,-1);
		\draw (e) -- (1.5,-1);
		\draw (f) -- (2.8,-1);
		\draw (f) -- (3.0,-1);
		\draw (f) -- (3.2,-1);
		\draw[black][dashed] (0,1) rectangle (4,-1.2);
	\end{tikzpicture}
	\caption{The support vertex $v$ is adjacent to $k\geq3$ leaves.}
	\label{lem3-fig1} 
\end{figure}
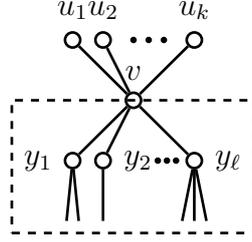

\pf Suppose that $N(v)=\{u_1,\cdots,u_k, y_1,\cdots, y_{\ell}\}$, where $d(u_i)=1$ for each $i\in\{1,\cdots,k\}$ and $d(y_j)\geq2$ for each $j\in \{1,\cdots,\ell\}$. See Figure \ref{lem3-fig1}. 
Let $T'=T-u_k.$

Note that $k\geq3$. We have
\begin{align*}
	\Phi(T)&=\sum\limits_{i=1}^k\Phi_{\{v,u_i\}}(T)+\sum\limits_{j=1}^{\ell}\Phi_{\{v,y_j\}}(T)+\Phi_{\overline{v}}(T),\\
	\Phi(T')=&\sum\limits_{i=1}^{k-1}\Phi_{\{v,u_i\}}(T')+\sum\limits_{j=1}^{\ell}\Phi_{\{v,y_j\}}(T')+\Phi_{\overline{v}}(T')\\
	=&\sum\limits_{i=1}^{k-1}\Phi_{\{v,u_i\}}(T)+\sum\limits_{j=1}^{\ell}\Phi_{\{v,y_j\}}(T)+\Phi_{\overline{v}}(T).
\end{align*}
It follows that
\[\Phi(T)=\Phi(T')+\Phi_{\{v,u_k\}}(T).\]

By Theorem \ref{thm1} and the fact that $k\geq3$, 
\[\Phi_{\{v,u_k\}}(T)\leq f_1(n-5)\leq 3^{\frac{n-5}{3}}.\]
Hence, if $\Phi(T')\leq t(n-1)$, then
\[\Phi(T)=\Phi(T')+\Phi_{\{v,u_k\}}(T)\leq t(n-1)+3^{\frac{n-5}{3}}<t(n).\]
\qed

Let $T^*_9$ be a tree of order $9$ obtained from $2K_{1,3}$ and an isolated vertex by adding two edges between the isolated vertex and two leaves of $2K_{1,3}$. By direct calculations, 
\[\Phi(T^*_7)=11=t(7),\ \Phi(T^*_8)=15<t(8), \ \Phi(T^*_9)=18<t(9).\]
\begin{lem}\label{lem4}
	For any integer $n\in\{7,8,9\}$ and any tree $T$ of order $n$, \[\Phi(T)\leq \Phi(T^*_n),\] with equality if and only if $T\cong T^*_n$.	
\end{lem}

\pf We can check readily that for any tree $T$ of order $6$, $\Phi(T)<t(6)$. For $n\in\{7,8,9\}$, by Lemma \ref{lem1}, \ref{lem2}, and \ref{lem3} and the fact that for any tree $T$ of order $n-1$, $\Phi(T)\leq t(n-1)$, we can verify the validity of Lemma \ref{lem4} for $n=7$, $8$, and $9$ successively. \qed

\section{Proof of Theorem \ref{thm2}}\label{sec3}

We prove both statements of Theorem \ref{thm2} simultaneously by induction on the order $n$ of trees and forests. 

{\bf Base cases.}

When $4\leq n\leq 6$, it is easy to prove that the results in both (1) and (2) hold true. 

(a) By Lemma \ref{lem4}, the results in (1) hold true for any integer $n\in \{7,8,9\}$. When $n=10$, let $T$ be a tree that has the most number of maximal dissociation sets among all trees of order $10$. We will prove that $\Phi(T)\leq t(10)$, with equality if and only if $T\cong T^*_{10}$. 

\begin{figure}[H]
	\begin{tikzpicture}
		[line width = 1pt, scale=0.8,
		empty/.style = {circle, draw, fill = white, inner sep=0mm, minimum size=2mm}, full/.style = {circle, draw, fill = black, inner sep=0mm, minimum size=2mm}, full1/.style = {circle, draw, fill = black, inner sep=0mm, minimum size=0.5mm}]
		\node [empty] (a) at (3.5,12.5) {};
		\node [empty] (b) at (4.5,12.5) {};
		\node [empty] (c) at (4,12) {};
		\node [empty] (d) at (4,11.5) {};
		\node [empty] (e) at (4.5,11.5) {};
		\node [empty] (f) at (5,11.75) {};
		\node [empty] (g) at (5,11.25) {};
		\node [empty] (h) at (4,11) {};
		\node [empty] (i) at (3.5,10.5) {};
		\node [empty] (j) at (4.5,10.5) {};
		\draw (a) -- (c);
		\draw (b) -- (c);
		\draw (c) -- (d);
		\draw (e) -- (d);
		\draw (e) -- (f);
		\draw (e) -- (g);
		\draw (d) -- (h);
		\draw (i) -- (h);
		\draw (h) -- (j);
		%\path (7,10) node[below]{$T_{10}(1)$};
		
		\node [empty] (a) at (6.5,13.5) {};
		\node [empty] (b) at (7.5,13.5) {};
		\node [empty] (c) at (7,13) {};
		\node [empty] (d) at (7,12.5) {};
		\node [empty] (e) at (7,12) {};
		\node [empty] (f) at (7.5,12) {};
		\node [empty] (g) at (7,11.5) {};
		\node [empty] (h) at (7,11) {};
		\node [empty] (i) at (6.5,10.5) {};
		\node [empty] (j) at (7.5,10.5) {};
		\draw (a) -- (c);
		\draw (b) -- (c);
		\draw (c) -- (d);
		\draw (e) -- (d);
		\draw (e) -- (f);
		\draw (e) -- (g);
		\draw (g) -- (h);
		\draw (i) -- (h);
		\draw (h) -- (j);
		%\path (7,10) node[below]{$T_{10}(1)$};
		
		\node [empty] (a) at (9.5,14) {};
		\node [empty] (b) at (10.5,14) {};
		\node [empty] (c) at (10,13.5) {};
		\node [empty] (d) at (10,13) {};
		\node [empty] (e) at (10,12.5) {};
		\node [empty] (f) at (10,12) {};
		\node [empty] (g) at (10,11.5) {};
		\node [empty] (h) at (10,11) {};
		\node [empty] (i) at (9.5,10.5) {};
		\node [empty] (j) at (10.5,10.5) {};
		\draw (a) -- (c);
		\draw (b) -- (c);
		\draw (c) -- (d);
		\draw (e) -- (d);
		\draw (e) -- (f);
		\draw (f) -- (g);
		\draw (g) -- (h);
		\draw (i) -- (h);
		\draw (h) -- (j);
		%\path (10,10) node[below]{$T_{10}(2)$};
		
	\end{tikzpicture}
	\caption{Three possible structures of $T$.}
	\label{thm2-fig1}
\end{figure}
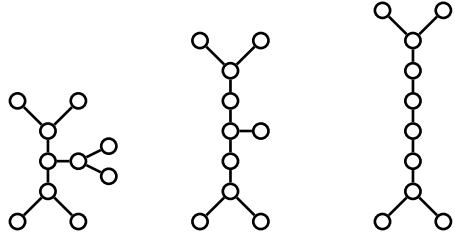

By Lemma \ref{lem3}, we can assume that any vertex of $T$ is adjacent to at most two leaves. By Lemma \ref{lem1} and \ref{lem2}, there are only three possible structures of $T$ which are pictured in Figure \ref{thm2-fig1}.

By direct calculations, 	
$\Phi(T)\leq t(10)$
with equality if and only if $T\cong T^*_{10}$.
\\

(b) For any integer $n\in\{7,8,9,10\}$, let $F$ be a forest that has the second largest number of maximal dissociation sets among all forests of order $n$.

When $n=7$, by Lemma \ref{lem4} and Theorem \ref{thm1}, it is easy to check that $\Phi(F)\leq f_2(n)$, with equality if and only if $F\cong T^*_7$.

When $n=8$, by Lemma \ref{lem4} and Theorem \ref{thm1}, it is easy to check that $\Phi(F)\leq f_2(n)$, with equality if and only if $F\cong T^*_8$ or $F\cong P_3\cup K_{1,4}$.

When $n=9$, because $\Phi(T^*_9)< \Phi(K_{1,3}\cup K_{1,4})$, there are at least two components in $F$. It is easy to check that $\Phi(F)\leq f_2(n)$, with equality if and only if $F\cong K_{1,3}\cup K_{1,4}$.

When $n=10$, because $t(10)<\Phi(P_3\cup T^*_7)$, there are at least two components in $F$. It is easy to check that $\Phi(F)\leq f_2(n)$, with equality if and only if $F\cong P_3\cup T^*_7$. 
\\
 
{\bf Inductive proof.}

Suppose that $n\geq11$ and for trees and forests of order at most $n-1$, the results in both (1) and (2) hold true.

(a) We first prove that the results in (1) hold true for trees of order $n$. 

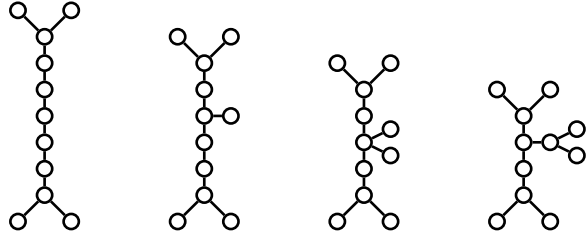
\begin{figure}[H]
	\begin{tikzpicture}
		[line width = 1pt, scale=0.7,
		empty/.style = {circle, draw, fill = white, inner sep=0mm, minimum size=2mm}, full/.style = {circle, draw, fill = black, inner sep=0mm, minimum size=2mm}, full1/.style = {circle, draw, fill = black, inner sep=0mm, minimum size=0.5mm}]
		\node [empty] (a) at (0.5,9) {};
		\node [empty] (b) at (1.5,9) {};
		\node [empty] (c) at (1,8.5) {};
		\node [empty] (d) at (1,8) {};
		\node [empty] (e) at (1,7.5) {};
		\node [empty] (f) at (1,7) {};
		\node [empty] (g) at (1,6.5) {};
		\node [empty] (h) at (1,6) {};
		\node [empty] (i) at (1,5.5) {};
		\node [empty] (j) at (0.5,5) {};
		\node [empty] (k) at (1.5,5) {};
		\draw (a) -- (c);
		\draw (b) -- (c);
		\draw (c) -- (d);
		\draw (e) -- (d);
		\draw (e) -- (f);
		\draw (f) -- (g);
		\draw (g) -- (h);
		\draw (i) -- (h);
		\draw (i) -- (k);
		\draw (i) -- (j);
		%\path (1,4.5) node[below]{$T_{11}(1)$};
		
		\node [empty] (a) at (3.5,8.5) {};
		\node [empty] (b) at (4.5,8.5) {};
		\node [empty] (c) at (4,8) {};
		\node [empty] (d) at (4,7.5) {};
		\node [empty] (e) at (4,7) {};
		\node [empty] (f) at (4.5,7) {};
		\node [empty] (g) at (4,6.5) {};
		\node [empty] (h) at (4,6) {};
		\node [empty] (i) at (4,5.5) {};
		\node [empty] (j) at (3.5,5) {};
		\node [empty] (k) at (4.5,5) {};
		\draw (a) -- (c);
		\draw (b) -- (c);
		\draw (c) -- (d);
		\draw (e) -- (d);
		\draw (e) -- (f);
		\draw (e) -- (g);
		\draw (g) -- (h);
		\draw (i) -- (h);
		\draw (i) -- (k);
		\draw (i) -- (j);
		%\path (4,4.5) node[below]{$T_{11}(2)$};
		
		\node [empty] (a) at (6.5,8) {};
		\node [empty] (b) at (7.5,8) {};
		\node [empty] (c) at (7,7.5) {};
		\node [empty] (d) at (7,7) {};
		\node [empty] (e) at (7,6.5) {};
		\node [empty] (f) at (7.5,6.75) {};
		\node [empty] (g) at (7.5,6.25) {};
		\node [empty] (h) at (7,6) {};
		\node [empty] (i) at (7,5.5) {};
		\node [empty] (j) at (6.5,5) {};
		\node [empty] (k) at (7.5,5) {};
		\draw (a) -- (c);
		\draw (b) -- (c);
		\draw (c) -- (d);
		\draw (e) -- (d);
		\draw (e) -- (f);
		\draw (e) -- (g);
		\draw (e) -- (h);
		\draw (i) -- (h);
		\draw (i) -- (k);
		\draw (i) -- (j);
		%\path (7,4.5) node[below]{$T_{11}(3)$};
		
		\node [empty] (a) at (9.5,7.5) {};
		\node [empty] (b) at (10.5,7.5) {};
		\node [empty] (c) at (10,7) {};
		\node [empty] (d) at (10,6.5) {};
		\node [empty] (e) at (10.5,6.5) {};
		\node [empty] (f) at (11,6.75) {};
		\node [empty] (g) at (11,6.25) {};
		\node [empty] (h) at (10,6) {};
		\node [empty] (i) at (10,5.5) {};
		\node [empty] (j) at (9.5,5) {};
		\node [empty] (k) at (10.5,5) {};
		\draw (a) -- (c);
		\draw (b) -- (c);
		\draw (c) -- (d);
		\draw (e) -- (d);
		\draw (e) -- (f);
		\draw (e) -- (g);
		\draw (d) -- (h);
		\draw (i) -- (h);
		\draw (i) -- (k);
		\draw (i) -- (j);
		%\path (7,4.5) node[below]{$T_{11}(3)$};
	\end{tikzpicture}
	\caption{Four possible structures of $T$.}
	\label{case1-fig1}
\end{figure}

When $n=11$, by Lemma \ref{lem1} and \ref{lem2}, there are only four possible structures of $T$ which are pictured in Figure \ref{case1-fig1}.

By direct calculations, 	
$\Phi(T)<t(n).$

\begin{figure}[H]
	\begin{tikzpicture}
		[line width = 1pt, scale=0.7,
		empty/.style = {circle, draw, fill = white, inner sep=0mm, minimum size=2mm}, full/.style = {circle, draw, fill = black, inner sep=0mm, minimum size=2mm}, full1/.style = {circle, draw, fill = black, inner sep=0mm, minimum size=0.5mm}]
		\node [empty] (a) at (0.5,9) {};
		\node [empty] (b) at (1.5,9) {};
		\node [empty] (c) at (1,8.5) {};
		\node [empty] (d) at (1,8) {};
		\node [empty] (e) at (1,7.5) {};
		\node [empty] (f) at (1,7) {};
		\node [empty] (g) at (1,6.5) {};
		\node [empty] (h) at (1,6) {};
		\node [empty] (i) at (1,5.5) {};
		\node [empty] (j) at (1,5) {};
		\node [empty] (k) at (0.5,4.5) {};
		\node [empty] (l) at (1.5,4.5) {};
		\draw (a) -- (c);
		\draw (b) -- (c);
		\draw (c) -- (d);
		\draw (e) -- (d);
		\draw (e) -- (f);
		\draw (f) -- (g);
		\draw (g) -- (h);
		\draw (i) -- (h);
		\draw (i) -- (j);
		\draw (j) -- (k);
		\draw (l) -- (j);
		%\path (1,4) node[below]{$T_{12}(1)$};
		
		\node [empty] (a) at (3.5,8.5) {};
		\node [empty] (b) at (4.5,8.5) {};
		\node [empty] (c) at (4,8) {};
		\node [empty] (d) at (4,7.5) {};
		\node [empty] (e) at (4,7) {};
		\node [empty] (f) at (4.5,7) {};
		\node [empty] (g) at (4,6.5) {};
		\node [empty] (h) at (4,6) {};
		\node [empty] (i) at (4,5.5) {};
		\node [empty] (j) at (4,5) {};
		\node [empty] (k) at (3.5,4.5) {};
		\node [empty] (l) at (4.5,4.5) {};
		\draw (a) -- (c);
		\draw (b) -- (c);
		\draw (c) -- (d);
		\draw (e) -- (d);
		\draw (e) -- (f);
		\draw (e) -- (g);
		\draw (g) -- (h);
		\draw (i) -- (h);
		\draw (i) -- (j);
		\draw (j) -- (k);
		\draw (l) -- (j);
		%\path (4,4) node[below]{$T_{12}(2)$};
		
		\node [empty] (a) at (6.5,8.5) {};
		\node [empty] (b) at (7.5,8.5) {};
		\node [empty] (c) at (7,8) {};
		\node [empty] (d) at (7,7.5) {};
		\node [empty] (e) at (7,7) {};
		\node [empty] (f) at (7,6.5) {};
		\node [empty] (g) at (7.5,6.5) {};
		\node [empty] (h) at (7,6) {};
		\node [empty] (i) at (7,5.5) {};
		\node [empty] (j) at (7,5) {};
		\node [empty] (k) at (6.5,4.5) {};
		\node [empty] (l) at (7.5,4.5) {};
		\draw (a) -- (c);
		\draw (b) -- (c);
		\draw (c) -- (d);
		\draw (e) -- (d);
		\draw (e) -- (f);
		\draw (f) -- (g);
		\draw (f) -- (h);
		\draw (i) -- (h);
		\draw (i) -- (j);
		\draw (j) -- (k);
		\draw (l) -- (j);
		%\path (7,4) node[below]{$T_{12}(3)$};
		
		\node [empty] (a) at (0.5,3) {};
		\node [empty] (b) at (1.5,3) {};
		\node [empty] (c) at (1,2.5) {};
		\node [empty] (d) at (1,2) {};
		\node [empty] (e) at (1,1.5) {};
		\node [empty] (f) at (1,1) {};
		\node [empty] (g) at (1.5,1.25) {};
		\node [empty] (h) at (1.5,0.75) {};
		\node [empty] (i) at (1,0.5) {};
		\node [empty] (j) at (1,0) {};
		\node [empty] (k) at (0.5,-0.5) {};
		\node [empty] (l) at (1.5,-0.5) {};
		\draw (a) -- (c);
		\draw (b) -- (c);
		\draw (c) -- (d);
		\draw (e) -- (d);
		\draw (e) -- (f);
		\draw (f) -- (g);
		\draw (f) -- (h);
		\draw (i) -- (f);
		\draw (i) -- (j);
		\draw (j) -- (k);
		\draw (l) -- (j);
		%\path (1,-1) node[below]{$T_{12}(4)$};
		
		\node [empty] (a) at (3.5,3) {};
		\node [empty] (b) at (4.5,3) {};
		\node [empty] (c) at (4,2.5) {};
		\node [empty] (d) at (4,2) {};
		\node [empty] (e) at (4,1.5) {};
		\node [empty] (f) at (4.5,1.5) {};
		\node [empty] (g) at (4.5,1) {};
		\node [empty] (h) at (4,1) {};
		\node [empty] (i) at (4,0.5) {};
		\node [empty] (j) at (4,0) {};
		\node [empty] (k) at (3.5,-0.5) {};
		\node [empty] (l) at (4.5,-0.5) {};
		\draw (a) -- (c);
		\draw (b) -- (c);
		\draw (c) -- (d);
		\draw (e) -- (d);
		\draw (e) -- (f);
		\draw (e) -- (h);
		\draw (g) -- (h);
		\draw (i) -- (h);
		\draw (i) -- (j);
		\draw (j) -- (k);
		\draw (l) -- (j);
		%\path (4,-1) node[below]{$T_{12}(5)$};
		
		\node [empty] (a) at (6.5,2.5) {};
		\node [empty] (b) at (7.5,2.5) {};
		\node [empty] (c) at (7,2) {};
		\node [empty] (d) at (7,1.5) {};
		\node [empty] (e) at (8,1.25) {};
		\node [empty] (f) at (8,0.75) {};
		\node [empty] (g) at (7.5,1) {};
		\node [empty] (h) at (7,1) {};
		\node [empty] (i) at (7,0.5) {};
		\node [empty] (j) at (7,0) {};
		\node [empty] (k) at (6.5,-0.5) {};
		\node [empty] (l) at (7.5,-0.5) {};
		\draw (a) -- (c);
		\draw (b) -- (c);
		\draw (c) -- (d);
		\draw (h) -- (d);
		\draw (g) -- (e);
		\draw (f) -- (g);
		\draw (g) -- (h);
		\draw (i) -- (h);
		\draw (i) -- (j);
		\draw (j) -- (k);
		\draw (l) -- (j);
		%\path (7,-1) node[below]{$T_{12}(6)$};
		
		\node [empty] (a) at (9.5,2.5) {};
		\node [empty] (b) at (10.5,2.5) {};
		\node [empty] (c) at (10,2) {};
		\node [empty] (d) at (10,1.5) {};
		\node [empty] (e) at (11,1.75) {};
		\node [empty] (f) at (11,1.25) {};
		\node [empty] (g) at (10.5,1.5) {};
		\node [empty] (h) at (10,1) {};
		\node [empty] (i) at (10,0.5) {};
		\node [empty] (j) at (10,0) {};
		\node [empty] (k) at (9.5,-0.5) {};
		\node [empty] (l) at (10.5,-0.5) {};
		\draw (a) -- (c);
		\draw (b) -- (c);
		\draw (c) -- (d);
		\draw (h) -- (d);
		\draw (g) -- (e);
		\draw (f) -- (g);
		\draw (g) -- (d);
		\draw (i) -- (h);
		\draw (i) -- (j);
		\draw (j) -- (k);
		\draw (l) -- (j);
		%\path (7,-1) node[below]{$T_{12}(6)$};
	\end{tikzpicture}
	\caption{Seven possible structures of $T$.}
	\label{case2-fig1}
\end{figure}
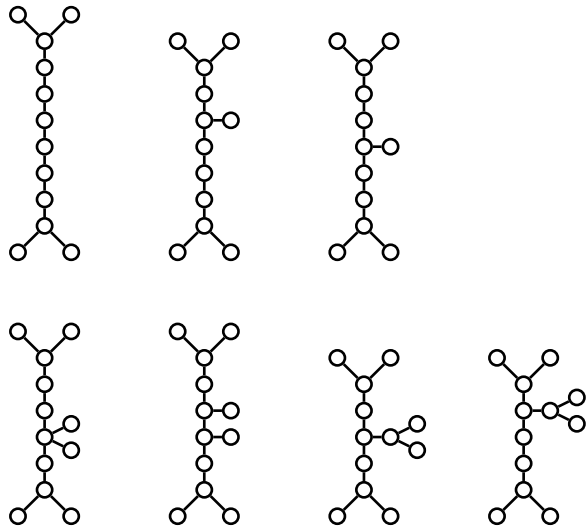

When $n=12$, by Lemma \ref{lem1} and \ref{lem2}, there are only seven possible structures of $T$ which are pictured in Figure \ref{case2-fig1}.

By direct calculations, 	
$\Phi(T)<t(n).$
\\

Now, we suppose that $n\geq13$. Let $T$ be a tree that has the most number of maximal dissociation sets among all trees of order $n$. We will prove that $\Phi(T)\leq t(n)$, with equality if and only if $n\equiv1\pmod{3}$ and $T\cong T^*_{n}$. 

If the diameter of $T$ is at most $3$, then $T$ has at most $n$ maximal dissociation sets. Consequently, we proceed under the assumption that the diameter of $T$ is at least 4.
By Lemma \ref{lem3} and the induction hypothesis, we assume that no vertex in $T$ is adjacent to more than two leaves.

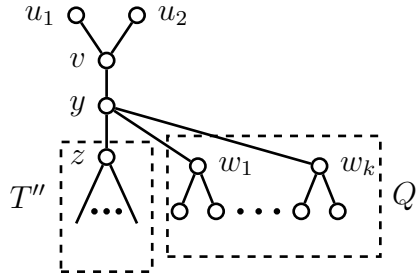
\begin{figure}[H]
				\begin{tikzpicture}
					[line width = 1pt, scale=0.8,
					empty/.style = {circle, draw, fill = white, inner sep=0mm, minimum size=2mm}, full/.style = {circle, draw, fill = black, inner sep=0mm, minimum size=2mm}, full1/.style = {circle, draw, fill = black, inner sep=0mm, minimum size=0.5mm}]
					\node [empty, label = left:$u_1$] (a) at (1.5,1.5) {};
					\node [empty, label = right:$u_2$] (b) at (2.5,1.5) {};
					\node [empty, label = left:$v$] (c) at (2,0.75) {};
					\node [empty, label = left:$y$] (d) at (2,0) {};
					%\node [empty] (e) at (1.5,-1.75) {};
					\node [full1] at (1.8,-1.75) {};
					\node [full1] at (2.0,-1.75) {};
					\node [full1] at (2.2,-1.75) {};
					\node [empty, label = left:$z$] (f) at (2,-0.85) {};
					%\node [empty] (g) at (2.5,-1.75) {};
					\node [empty, label = right:$w_1$] (h) at (3.5,-1) {};
					\node [empty] (i) at (3.2,-1.75) {};
					\node [empty] (j) at (3.8,-1.75) {};
					\node [full1] at (4.2,-1.75) {};
					\node [full1] at (4.5,-1.75) {};
					\node [full1] at (4.8,-1.75) {};
					\node [empty, label = right:$w_k$] (k) at (5.5,-1) {};
					\node [empty] (l) at (5.2,-1.75) {};
					\node [empty] (m) at (5.8,-1.75) {};
					\draw (a) -- (c);
					\draw (b) -- (c);
					\draw (c) -- (d);
					\draw (d) -- (f);
					\draw (1.5,-1.95) -- (f);
					\draw (f) -- (2.5,-1.95);
					\draw (d) -- (h);
					\draw (h) -- (i);
					\draw (h) -- (j);
					\draw (d) -- (k);
					\draw (k) -- (l);
					\draw (k) -- (m);
					\draw[black][dashed] (1.25,-0.6) rectangle (2.75,-2.75);
					\draw[black][dashed] (3,-0.5) rectangle (6.5,-2.5);
					\path (1.2,-1.5) node[left]{$T''$};
					\path (6.5,-1.5) node[right]{$Q$};
				\end{tikzpicture}
				\caption{there is at least two components in $T-\{u_1,u_2,v,y\}$.}
				\label{thm2-fig2}
			\end{figure}

Let $P:=u_1vyz\cdots$ be a longest path in $T$, where $u_1$ is a leaf. According to Lemma \ref{lem1}, $d(v)\neq2$, indicating that $v$ is adjacent to exactly two leaves. We designate the other leaf adjacent to $v$ as $u_2$. Since $P$ is a longest path in $T$, the vertex $y$ is the unique non-leaf vertex that is connected to $v$. Furthermore, Lemma \ref{lem2} asserts that the vertex $y$ cannot be a support vertex. See Figure \ref{thm2-fig2}.

By Lemma \ref{lem1} and \ref{lem2} and the fact that $P$ is a longest path in $T$, there is at most one component in $T-\{u_1,u_2,v,y\}$ that is not isomorphic to $P_3$. If any component in $T-\{u_1,u_2,v,y\}$ is isomorphic to $P_3$, then $n\equiv 1\pmod3$, $T\cong T^*_n$, and 
\[\Phi(T)=t(n).\]

Now, we proceed under the assumption that there exists exactly one component in $T-\{u_1,u_2,v,y\}$ that is not isomorphic to $P_3$.
It is easy to see that the vertex $z$ must be in this unique component.

\begin{cl}\label{cl1}
If there are at least two components in $T-\{u_1,u_2,v,y\}$ (see Figure \ref{thm2-fig2}), then \[\Phi(T)<t(n).\]
\end{cl}

\pf Suppose that $N(y)=\{v,z,w_1,\cdots,w_k\}$, where $k\geq1$ and for each $i\in\{1,\cdots,k\}$, $w_i$ is the non-leaf vertex of a $P_3$. See Figure \ref{thm2-fig2}. 

Let $Q$ denote the set of vertices containing $w_1,\cdots,w_k$ along with all $2k$ leaves adjacent to at least one $w_i$. Let $T'=T-Q$ and $T''=T'-\{u_1,u_2,v,y\}$. Note that $T''\ncong P_3$.

It is easy to see that 
\[\big|\{S\in MD(T): S\cap\{w_1,\cdots,w_k\}=\emptyset\}\big|=\Phi(T'),\]
\[\big|\{S\in MD(T): S\cap\{w_1,\cdots,w_k\}\neq\emptyset\text{\ and\ }y\notin S \}\big|=3(3^k-1)\Phi(T''),\]
and
\[\big|\{S\in MD(T): \big|S\cap\{w_1,\cdots,w_k\}\big|=1\text{\ and\ }y\in S \}\big|=k\Phi(T''-z).\]
Thus,
\[\Phi(T)=\Phi(T')+3(3^k-1)\Phi(T'')+k\Phi(T''-{z}).\]
		
Since $P=u_1vyz\cdots$ is a longest path in $T$, $V(T'')\geq4$. For the convenience of calculations, we suppose that $|V(T'')|=3p+1$, where $3p\geq3$ is an integer and $p$ may not be an integer. Thus,
\[n=3k+3p+5\geq13.\]

%Let $T_{(2)}$ have $3t+1$ vertices. We can see $t\geq\frac{2}{3}$ and $t=\frac{1}{3}x$($x$ is an integer greater than or equal to 2).

By the induction hypothesis and Theorem \ref{thm1}, we have
\begin{align*}
	\Phi(T)&=\Phi(T')+3(3^k-1)\Phi(T'')+k\Phi(T''-{z})\\
	&\leq t(3p+5)+3(3^k-1)t(3p+1)+kf_1(3p)\\
	&\leq 3^{p+\frac{4}{3}}+p+\frac{4}{3}+3(3^k-1)(3^{p}+p)+k3^p.
\end{align*}
In order to prove that 
			\[\Phi(T)< t(n)=3^{p+k+\frac{4}{3}}+p+k+\frac{4}{3},\]
we consider the following function
\begin{align*}
h(p,k):&=3^{p+k+\frac{4}{3}}+p+k+\frac{4}{3}-\left[(3^{p+\frac{4}{3}}+p+\frac{4}{3})+3(3^k-1)(3^{p}+p)+k 3^p\right]\\
&=\left[(3^{\frac{1}{3}}-1-\frac{p}{3^p})(3^{k+1}-3)-k\right]\cdot 3^p+k\\
&=\left[(3^{\frac{1}{3}}-1-\frac{p}{3^p})-\frac{k}{3^{k+1}-3}\right]\cdot (3^{k+1}-3) \cdot 3^p+k.
\end{align*}
Let
			\[g(p,k)=(3^{\frac{1}{3}}-1-\frac{p}{3^p})-\frac{k}{3^{k+1}-3}.\]

When $k=1$, $p\geq\frac{5}{3}$ and the function $\frac{p}{3^p}$ is monotonically decreasing. Hence, 
\begin{equation*}%\label{eqn1}
	g(p,k)=(3^{\frac{1}{3}}-1-\frac{p}{3^p})-\frac{k}{3^{k+1}-3}\geq (3^{\frac{1}{3}}-1-\frac{4}{3^{4}})-\frac{1}{6} > 0,
\end{equation*}
which implies $h(p,k)>0$. 
			
When $k\geq 2$, $\frac{k}{3^{k+1}-3}\leq \frac{1}{12}$ and for $p\geq1$, $\frac{p}{3^p}\leq \frac{1}{3}$.
So
			\[g(p,k)=(3^{\frac{1}{3}}-1-\frac{p}{3^p})-\frac{k}{3^{k+1}-3}\geq (3^{\frac{1}{3}}-1-\frac{1}{3})-\frac{1}{12} > 0,\]
which implies $h(p,k)>0$.
			
Thus, we have proved that $\Phi(T)<t(n)$. The proof of Claim \ref{cl1} is complete. \qed

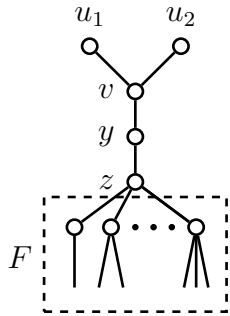
\begin{figure}[H]
				\begin{tikzpicture}
					[line width = 1pt, scale=0.8,
					empty/.style = {circle, draw, fill = white, inner sep=0mm, minimum size=2mm}, full/.style = {circle, draw, fill = black, inner sep=0mm, minimum size=2mm}, full1/.style = {circle, draw, fill = black, inner sep=0mm, minimum size=0.5mm}]
					\node [empty,label=above:$u_1$] (a) at (0.25,13) {};
					\node [empty,label=above:$u_2$] (b) at (1.75,13) {};
					\node [empty,label=left:$v$] (c) at (1,12.25) {};
					\node [empty,label=left:$y$] (d) at (1,11.5) {};
					\node [empty,label=left:$z$] (e) at (1,10.75) {};
					\node [empty] (f) at (0.6,10) {};
					\node [empty] (g) at (0,10) {};
					\node [full1] at (1,10) {};
					\node [full1] at (1.3,10) {};
					\node [full1] at (1.6,10) {};
					\node [empty] (h) at (2,10) {};
					\draw (a) -- (c);
					\draw (b) -- (c);
					\draw (c) -- (d);
					\draw (e) -- (d);
					\draw (e) -- (f);
					\draw (e) -- (g);
					\draw (e) -- (h);
					\draw (e) -- (h);
					\draw (f) -- (0.4,9);
					\draw (f) -- (0.8,9);
					\draw (g) -- (0,9);
					\draw (h) -- (1.8,9);
					\draw (h) -- (2.0,9);
					\draw (h) -- (2.2,9);
					\draw[black][dashed] (-0.5,10.5) rectangle (2.5,8.6);
					\path (-0.5,9.5) node[left]{$F$};
				\end{tikzpicture}
				\caption{There is exactly one component in $T-\{u_1,u_2,v,y\}$.}
				\label{thm2-fig3}
			\end{figure}

\begin{cl}\label{cl2}
	If there is exactly one component in $T-\{u_1,u_2,v,y\}$ (see Figure \ref{thm2-fig3}), then \[\Phi(T)<t(n).\]
\end{cl}

\pf The proof for the case when $13\leq n\leq19$ is quite cumbersome and has been included in Appendix for the sake of brevity. Now, we assume that $n\geq 20$.

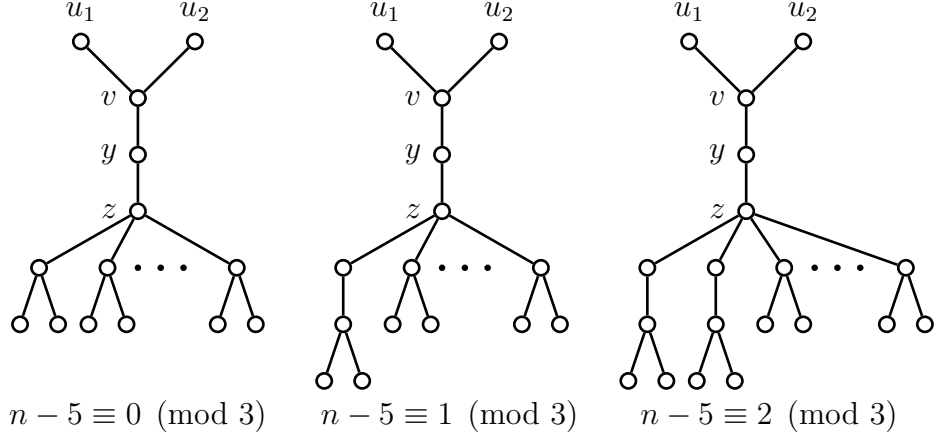
\begin{figure}[H]
	\begin{tikzpicture}
		[line width = 1pt, 
		empty/.style = {circle, draw, fill = white, inner sep=0mm, minimum size=2mm}, full/.style = {circle, draw, fill = black, inner sep=0mm, minimum size=2mm}, full1/.style = {circle, draw, fill = black, inner sep=0mm, minimum size=0.5mm}]
		\node [empty,label=above:$u_1$] (a) at (0.25,13) {};
		\node [empty,label=above:$u_2$] (b) at (1.75,13) {};
		\node [empty,label=left:$v$] (c) at (1,12.25) {};
		\node [empty,label=left:$y$] (d) at (1,11.5) {};
		\node [empty,label=left:$z$] (e) at (1,10.75) {};
		\node [empty] (f) at (0.6,10) {};
		\node [empty] (g) at (-0.3,10) {};
		\node [full1] at (1,10) {};
		\node [full1] at (1.3,10) {};
		\node [full1] at (1.6,10) {};
		\node [empty] (h) at (2.3,10) {};
		\node [empty] (i) at (-0.55,9.25) {};
		\node [empty] (j) at (-0.05,9.25) {};
		\node [empty] (k) at (0.35,9.25) {};
		\node [empty] (l) at (0.85,9.25) {};
		\node [empty] (m) at (2.05,9.25) {};
		\node [empty] (n) at (2.55,9.25) {};
		\draw (a) -- (c);
		\draw (b) -- (c);
		\draw (c) -- (d);
		\draw (e) -- (d);
		\draw (e) -- (f);
		\draw (e) -- (g);
		\draw (e) -- (h);
		\draw (g) -- (i);
		\draw (g) -- (j);
		\draw (f) -- (k);
		\draw (f) -- (l);
		\draw (h) -- (m);
		\draw (h) -- (n);
		\path (1,8.4) node[below]{$n-5\equiv 0\pmod3$};
		\node [empty,label=above:$u_1$] (a) at (4.25,13) {};
		\node [empty,label=above:$u_2$] (b) at (5.75,13) {};
		\node [empty,label=left:$v$] (c) at (5,12.25) {};
		\node [empty,label=left:$y$] (d) at (5,11.5) {};
		\node [empty,label=left:$z$] (e) at (5,10.75) {};
		\node [empty] (f) at (4.6,10) {};
		\node [empty] (g) at (3.7,10) {};
		\node [full1] at (5,10) {};
		\node [full1] at (5.3,10) {};
		\node [full1] at (5.6,10) {};
		\node [empty] (h) at (6.3,10) {};
		\node [empty] (o) at (3.7,9.25) {};
		\node [empty] (i) at (3.45,8.5) {};
		\node [empty] (j) at (3.95,8.5) {};
		\node [empty] (k) at (4.35,9.25) {};
		\node [empty] (l) at (4.85,9.25) {};
		\node [empty] (m) at (6.05,9.25) {};
		\node [empty] (n) at (6.55,9.25) {};
		\draw (a) -- (c);
		\draw (b) -- (c);
		\draw (c) -- (d);
		\draw (e) -- (d);
		\draw (e) -- (f);
		\draw (e) -- (g);
		\draw (e) -- (h);
		\draw (g) -- (o);
		\draw (o) -- (i);
		\draw (o) -- (j);
		\draw (f) -- (k);
		\draw (f) -- (l);
		\draw (h) -- (m);
		\draw (h) -- (n);
		\path (5.1,8.4) node[below]{$n-5\equiv 1\pmod3$};
		\node [empty,label=above:$u_1$] (a) at (8.25,13) {};
		\node [empty,label=above:$u_2$] (b) at (9.75,13) {};
		\node [empty,label=left:$v$] (c) at (9,12.25) {};
		\node [empty,label=left:$y$] (d) at (9,11.5) {};
		\node [empty,label=left:$z$] (e) at (9,10.75) {};
		\node [empty] (f) at (8.6,10) {};
		\node [empty] (g) at (7.7,10) {};
		\node [full1] at (9.9,10) {};
		\node [full1] at (10.2,10) {};
		\node [full1] at (10.5,10) {};
		\node [empty] (h) at (9.5,10) {};
		\node [empty] (o) at (7.7,9.25) {};
		\node [empty] (i) at (7.45,8.5) {};
		\node [empty] (j) at (7.95,8.5) {};
		\node [empty] (p) at (8.6,9.25) {};
		\node [empty] (k) at (8.35,8.5) {};
		\node [empty] (l) at (8.85,8.5) {};
		\node [empty] (m) at (9.25,9.25) {};
		\node [empty] (n) at (9.75,9.25) {};
		\node [empty] (q) at (11.1,10) {};
		\node [empty] (r) at (10.85,9.25) {};
		\node [empty] (s) at (11.35,9.25) {};
		\draw (a) -- (c);
		\draw (b) -- (c);
		\draw (c) -- (d);
		\draw (e) -- (d);
		\draw (e) -- (f);
		\draw (e) -- (g);
		\draw (e) -- (h);
		\draw (g) -- (o);
		\draw (o) -- (i);
		\draw (o) -- (j);
		\draw (f) -- (p);
		\draw (p) -- (k);
		\draw (p) -- (l);
		\draw (h) -- (m);
		\draw (h) -- (n);
		\draw (e) -- (q);
		\draw (q) -- (r);
		\draw (q) -- (s);
		\path (9.3,8.4) node[below]{$n-5\equiv 2\pmod3$};
	\end{tikzpicture}
\caption{When $\Phi(F)=f_1(n)$, the structure of $T$ can be determined.}
\label{thm2-fig4}
\end{figure}

Let $F:=T-\{u_1,u_2,v,y,z\}$. So, $|V(F)|=n-5$. By Theorem \ref{thm1}, if $\Phi(F)=f_1(n-5)$, the structure of $F$ is uniquely determined. Moreover, by Lemma \ref{lem1} and \ref{lem2}, we can unambiguously determine the structure of $T$ as depicted in Figure \ref{thm2-fig4}.

By direct calculations, when $n\geq20$,

(1) if $n-5\equiv 0\pmod3$, then $\Phi(T)=4\cdot 3^{\frac{n-5}{3}}+n-5<t(n),$

(2) if $n-5\equiv 1\pmod3$, then $\Phi(T)=16\cdot 3^{\frac{n-9}{3}}+3n-25<t(n),$

(3) if $n-5\equiv 2\pmod3$, then $\Phi(T)=64\cdot 3^{\frac{n-13}{3}}+9n-82<t(n).$

Now, we turn our attention to the case when $\Phi(F)<f_1(n-5)$. Note that when $n\geq20$, $f_2(n)\leq 11\cdot 3^{\frac{n-7}{3}}$. By the induction hypothesis, $\Phi(F)\leq f_2(n-5)\leq 11\cdot 3^{\frac{n-12}{3}}$ and 
	\begin{align*}
	\Phi(T)=&\Phi_{\{u_1,u_2\}}(T)+2\cdot \Phi_{\{u_1,v\}}(T)+\Phi_{\{v,y\}}(T)\\
	\leq &t(n-3)+2\cdot t(n-4)+\Phi(F)\\
	\leq &3^{\frac{n-4}{3}}+\frac{n-4}{3}+2\cdot (3^{\frac{n-5}{3}}+\frac{n-5}{3})+11\cdot 3^{\frac{n-12}{3}}\\
	<&t(n).
\end{align*}

The proof of Claim \ref{cl2} is complete. \qed

Thus, by the induction hypothesis, we have proved that \[\Phi(T)\leq t(n),\]
with equality if and only if $n\equiv 1\pmod3$ and $T\cong T^*_n$.

\vskip 3mm

(b) Next, we prove that the results in (2) hold true for forests of order $n$. Suppose that $n\geq11$, the results in (2) hold true for forests of order at most $n-1$, and the results in $(1)$ hold true for trees of order at most $n$.

Let $F$ be a forest of order $n$ that has the most number of maximal dissociation sets among all forests of order $n$ satisfying $\Phi(F)<f_1(n)$. We will prove that $\Phi(F)\leq f_2(n)$.

If $F$ is connected, then it follows that $\Phi(F)\leq t(n)<f_2(n)$. Thus, we assume that there are at least two components in $F$. If $F$ contains a component of order at most $2$, then Theorem \ref{thm1} implies that
\[\Phi(F)\leq f_1(n-1)\leq3^{\frac{n-1}{3}}<f_2(n).\]
Hence, we can further assume that the order of any component in $F$ is at least 3.

\begin{cl}\label{cl4}
	If any of the following three conditions is satisfied:
	
	(1) $n\equiv 1\pmod3$ and there is a component in $F$ of order at least $8$,
	
	(2) $n\equiv 2\pmod3$ and there is a component in $F$ of order at least $9$,
	
	(3) $n\equiv 0\pmod3$ and there is a component in $F$ of order at least $10$,
	
\noindent then \[\Phi(F)< f_2(n).\]	
\end{cl}
\pf If $n\equiv 1\pmod3$ and there is a component in $F$ of order $n_1\geq8$, then by Theorem \ref{thm1} and Lemma \ref{lem4},
\[\Phi(F)\leq t(n_1)\cdot 3^{\frac{n-n_1}{3}}\leq15\cdot3^{\frac{n-8}{3}}<11\cdot 3^{\frac{n-7}{3}}.\]

If $n\equiv 2\pmod3$ and there is a component in $F$ of order $n_2\geq9$, then 
\[\Phi(F)\leq t(n_2)\cdot 3^{\frac{n-n_2}{3}}\leq 18\cdot3^{\frac{n-9}{3}}< 15\cdot 3^{\frac{n-8}{3}}.\]

If $n\equiv 0\pmod3$ and there is a component in $F$ of order $n_3\geq10$, then
\[\Phi(F)\leq t(n_3)\cdot 3^{\frac{n-n_3}{3}}\leq 30\cdot 3^{\frac{n-10}{3}} < 4^3\cdot 3^{\frac{n-12}{3}}.\]

Thus, if any of the three conditions is satisfied, then $\Phi(F)< f_2(n).$ \qed
\\

Let $r$ be the maximum order among all components in $F$. We consider the following three cases.

{\bf Case 1.}	$n\equiv 1\pmod3$.

By Claim \ref{cl4}, we can assume that $r\leq 7$. By Theorem \ref{thm1} and Lemma \ref{lem4}, when $r=7$, \[\Phi(F)\leq 11\cdot 3^{\frac{n-7}{3}},\]
with equality if and only if $F\cong T^*_7\cup \frac{n-7}{3}P_3$. When $r=6$,
\[\Phi(F)\leq 6\cdot 4\cdot 3^{\frac{n-10}{3}}< 11\cdot 3^{\frac{n-7}{3}}.\] 
When $r=5$,
\[\Phi(F)\leq 5\cdot 4^2\cdot 3^{\frac{n-13}{3}}< 11\cdot 3^{\frac{n-7}{3}}.\]
Finally, when $r=4$, because $\Phi(F)<f_1(n)$,  
	\[\Phi(F)\leq 4^4\cdot 3^{\frac{n-16}{3}}< 11\cdot 3^{\frac{n-7}{3}}.\]
			
In summary, for the case when $n\equiv 1\pmod3$, we have \[\Phi(F)\leq 11\cdot 3^{\frac{n-7}{3}},\]
with equality if and only if $F\cong T^*_7\cup \frac{n-7}{3}P_3$.
\\

{\bf Case 2.}	$n\equiv 2\pmod3$.

By Claim \ref{cl4}, we can assume that $r\leq 8$. By Theorem \ref{thm1} and Lemma \ref{lem4}, when $r=8$,
\[\Phi(F)\leq 15\cdot 3^{\frac{n-8}{3}},\]
with equality if and only if $F\cong T^*_8\cup \frac{n-8}{3}P_3$. When $r=7$,
	\[\Phi(F)\leq 11\cdot 4\cdot 3^{\frac{n-11}{3}}< 15\cdot 3^{\frac{n-8}{3}}.\]
When $r=6$,
\[\Phi(F)\leq 6\cdot 4^2\cdot 3^{\frac{n-14}{3}}< 15\cdot 3^{\frac{n-8}{3}}.\]
When $r=5$,
\[\Phi(F)\leq 5\cdot 3^{\frac{n-5}{3}}= 15\cdot 3^{\frac{n-8}{3}},\]
with equality if and only if $F\cong K_{1,4}\cup \frac{n-5}{3}P_3$.
Finally, when $r=4$, because $\Phi(F)<f_1(n)$,  
	\[\Phi(F)\leq 4^5\cdot 3^{\frac{n-20}{3}}< 15\cdot 3^{\frac{n-8}{3}}.\]

In summary, for the case when $n\equiv 2\pmod3$, we have \[\Phi(F)\leq 15\cdot 3^{\frac{n-8}{3}},\]
with equality if and only if $F\cong T^*_8\cup \frac{n-8}{3}P_3$ or $F\cong K_{1,4}\cup \frac{n-5}{3}P_3$.
\\

{\bf Case 3.} $n\equiv 0\pmod3$.
		
By Claim \ref{cl4} and the fact that $\Phi(F)<f_1(n)$, we can assume that $4\leq r\leq9$. By Theorem \ref{thm1} and Lemma \ref{lem4}, when $r=9$,
\[\Phi(F)\leq 18\cdot 3^{\frac{n-9}{3}}< 4^3\cdot 3^{\frac{n-12}{3}}.\]
When $r=8$,
\[\Phi(F)\leq 15\cdot 4\cdot 3^{\frac{n-12}{3}}< 4^3\cdot 3^{\frac{n-12}{3}}.\]
When $r=7$,
\[\Phi(F)\leq 11\cdot 4^2\cdot 3^{\frac{n-15}{3}}< 4^3\cdot 3^{\frac{n-12}{3}}.\]
When $r=6$,
\[\Phi(F)\leq 6\cdot 3^{\frac{n-6}{3}}< 4^3\cdot 3^{\frac{n-12}{3}}.\]
When $r=5$,		
\[\Phi(F)\leq 5\cdot 4\cdot 3^{\frac{n-9}{3}}< 4^3\cdot 3^{\frac{n-12}{3}}.\]
When $r=4$, 		
\[\Phi(F)\leq 4^3\cdot 3^{\frac{n-12}{3}}= 4^3\cdot 3^{\frac{n-12}{3}}.\]
with equality if and only if $F\cong 3K_{1,3}\cup \frac{n-12}{3}P_3$.
			
In summary, for the case when $n\equiv 0\pmod3$, we have \[\Phi(F)\leq 4^3\cdot 3^{\frac{n-12}{3}},\]
with equality if and only if $F\cong 3K_{1,3}\cup \frac{n-12}{3}P_3$.

Thus, we have proved that the results in (2) hold true for forests of order $n$.			

We complete the proof of Theorem \ref{thm2}. \qed

\section{Conclusion}\label{sec4}

In this paper, we present an upper bound for the number of maximal dissociation sets in trees of order $n$ and identify the tree that attains this bound. When $n\equiv1\pmod3$, the upper bound coincides with the largest number of maximal dissociation sets in trees of order $n$. But when $n\equiv 0\pmod3$ and $n\equiv 2\pmod3$, the largest number of maximal dissociation sets in trees of order $n$ remains unknown. We pose the following conjecture.

\begin{conj}\label{conj1}
The largest number of maximal dissociation sets in trees of order $n\geq7$ is 
\begin{align*}
	\begin{split}
		\left \{
		\begin{array}{ll}
			%5    & \hbox{if $n=5$,} \\
		%	6    & \hbox{if $n=6$,} \\
			3^{\frac{n-1}{3}}+\frac{n-1}{3}, & \hbox{if $n\equiv1\pmod{3}$,} \\
			4\cdot 3^{\frac{n-5}{3}}+n-5, &\hbox{if $n\equiv2\pmod3$,}\\
			16\cdot 3^{\frac{n-9}{3}}+3n-25, &\hbox{if $n\equiv0\pmod3$.}\\
		\end{array}
		\right.
	\end{split}
\end{align*}
The extremal trees that achieve this largest number are depicted in Figure \ref{conj-fig1}. Notably, when $n=14$, there exist two distinct extremal trees.
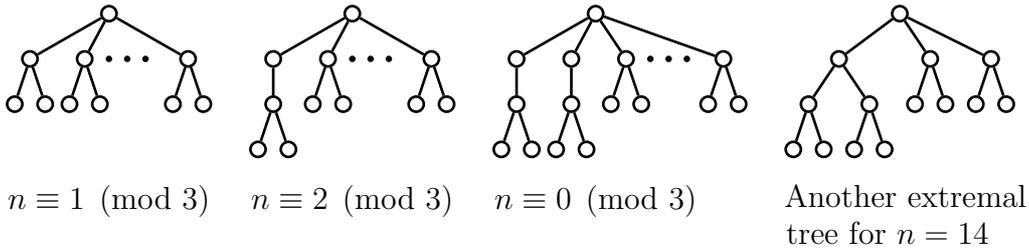
\begin{figure}[H]
	\begin{tikzpicture}
		[line width = 1pt, scale=0.8,
		empty/.style = {circle, draw, fill = white, inner sep=0mm, minimum size=2mm}, full/.style = {circle, draw, fill = black, inner sep=0mm, minimum size=2mm}, full1/.style = {circle, draw, fill = black, inner sep=0mm, minimum size=0.5mm}]
		
		\node [empty] (e) at (1,10.75) {};
		\node [empty] (f) at (0.6,10) {};
		\node [empty] (g) at (-0.3,10) {};
		\node [full1] at (1,10) {};
		\node [full1] at (1.3,10) {};
		\node [full1] at (1.6,10) {};
		\node [empty] (h) at (2.3,10) {};
		\node [empty] (i) at (-0.55,9.25) {};
		\node [empty] (j) at (-0.05,9.25) {};
		\node [empty] (k) at (0.35,9.25) {};
		\node [empty] (l) at (0.85,9.25) {};
		\node [empty] (m) at (2.05,9.25) {};
		\node [empty] (n) at (2.55,9.25) {};
		\draw (e) -- (f);
		\draw (e) -- (g);
		\draw (e) -- (h);
		\draw (g) -- (i);
		\draw (g) -- (j);
		\draw (f) -- (k);
		\draw (f) -- (l);
		\draw (h) -- (m);
		\draw (h) -- (n);
		\path (1,8.1) node[below]{$n\equiv1\pmod3$};
		
		\node [empty] (e) at (5,10.75) {};
		\node [empty] (f) at (4.6,10) {};
		\node [empty] (g) at (3.7,10) {};
		\node [full1] at (5,10) {};
		\node [full1] at (5.3,10) {};
		\node [full1] at (5.6,10) {};
		\node [empty] (h) at (6.3,10) {};
		\node [empty] (o) at (3.7,9.25) {};
		\node [empty] (i) at (3.45,8.5) {};
		\node [empty] (j) at (3.95,8.5) {};
		\node [empty] (k) at (4.35,9.25) {};
		\node [empty] (l) at (4.85,9.25) {};
		\node [empty] (m) at (6.05,9.25) {};
		\node [empty] (n) at (6.55,9.25) {};
		\draw (e) -- (f);
		\draw (e) -- (g);
		\draw (e) -- (h);
		\draw (g) -- (o);
		\draw (o) -- (i);
		\draw (o) -- (j);
		\draw (f) -- (k);
		\draw (f) -- (l);
		\draw (h) -- (m);
		\draw (h) -- (n);
		\path (5,8.1) node[below]{$n\equiv2\pmod3$};

		\node [empty] (e) at (9,10.75) {};
		\node [empty] (f) at (8.6,10) {};
		\node [empty] (g) at (7.7,10) {};
		\node [full1] at (9.9,10) {};
		\node [full1] at (10.2,10) {};
		\node [full1] at (10.5,10) {};
		\node [empty] (h) at (9.5,10) {};
		\node [empty] (o) at (7.7,9.25) {};
		\node [empty] (i) at (7.45,8.5) {};
		\node [empty] (j) at (7.95,8.5) {};
		\node [empty] (p) at (8.6,9.25) {};
		\node [empty] (k) at (8.35,8.5) {};
		\node [empty] (l) at (8.85,8.5) {};
		\node [empty] (m) at (9.25,9.25) {};
		\node [empty] (n) at (9.75,9.25) {};
		\node [empty] (q) at (11.1,10) {};
		\node [empty] (r) at (10.85,9.25) {};
		\node [empty] (s) at (11.35,9.25) {};
		\draw (e) -- (f);
		\draw (e) -- (g);
		\draw (e) -- (h);
		\draw (g) -- (o);
		\draw (o) -- (i);
		\draw (o) -- (j);
		\draw (f) -- (p);
		\draw (p) -- (k);
		\draw (p) -- (l);
		\draw (h) -- (m);
		\draw (h) -- (n);
		\draw (e) -- (q);
		\draw (q) -- (r);
		\draw (q) -- (s);
		\path (9,8.1) node[below]{$n\equiv0\pmod3$};

		\node [empty] (b) at (13,10) {};
		\node [empty] (c) at (14,10.75) {};
		\node [empty] (d) at (12.5,9.25) {};
		\node [empty] (e) at (13.5,9.25) {};
		\node [empty] (f) at (14.5,10) {};
		\node [empty] (g) at (15.5,10) {};
		\node [empty] (h) at (12.25,8.5) {};
		\node [empty] (i) at (12.75,8.5) {};
		\node [empty] (j) at (13.25,8.5) {};
		\node [empty] (k) at (13.75,8.5) {};
		\node [empty] (l) at (14.25,9.25) {};
		\node [empty] (m) at (14.75,9.25) {};
		\node [empty] (n) at (15.25,9.25) {};
		\node [empty] (o) at (15.75,9.25) {};
		\draw (b) -- (c);
		\draw (b) -- (d);
		\draw (b) -- (e);
		\draw (c) -- (f);
		\draw (c) -- (g);
		\draw (d) -- (h);
		\draw (d) -- (i);
		\draw (e) -- (j);
		\draw (e) -- (k);
		\draw (f) -- (l);
		\draw (f) -- (m);
		\draw (g) -- (n);
		\draw (g) -- (o);
		\path (14.1,8.1) node[below]{Another extremal};
		\path (13.8,7.5) node[below]{tree for $n=14$};
	\end{tikzpicture}
	\caption{The extremal trees.}
	\label{conj-fig1}
\end{figure}

\end{conj}

\section*{Appendix}

In the appendix we give the proof of Claim \ref{cl2} in the proof of Theorem \ref{thm2} for the case when $13\leq n\leq19$.

Now, there is exactly one component in $T-\{u_1,u_2,v,y\}$ and $d(y)=2$.

{\bf Case 1. $d(z)>2$.}

\begin{figure}[H]
		\begin{tikzpicture}
			[line width = 1pt, scale=0.8,
			empty/.style = {circle, draw, fill = white, inner sep=0mm, minimum size=2mm}, full/.style = {circle, draw, fill = black, inner sep=0mm, minimum size=2mm}, full1/.style = {circle, draw, fill = black, inner sep=0mm, minimum size=0.5mm}]
			\node [empty, label = above:$u_1$] (a) at (0.5,3) {};
			\node [empty, label=above:$u_2$] (b) at (1.5,3) {};
			\node [empty, label=left:$v$] (c) at (1,2.5) {};
			\node [empty, label=left:$y$] (d) at (1,2) {};
			\node [empty, label=left:$z$] (e) at (1,1.5) {};
			\node [empty] (f) at (1.5,2) {};
			\node [full1] at (0.8,0.6) {};
			\node [full1] at (1.0,0.6) {};
			\node [full1] at (1.2,0.6) {};
			%\node [empty] (i) at (1,1) {};
			%\node [empty] (j) at (0.5,1) {};
			%\node [empty] (k) at (1.5,1) {};
			\draw (a) -- (c);
			\draw (b) -- (c);
			\draw (c) -- (d);
			\draw (e) -- (d);
			\draw (e) -- (f);
			%\draw (e) -- (i);
			%\draw (e) -- (j);
			%\draw (e) -- (k);
			\draw (e) -- (0.5,0.5);
			%\draw (i) -- (1.1,0.25);
			%\draw (j) -- (0.5,0.25);
			\draw (e) -- (1.5,0.5);
			\draw[black][dashed] (0.25,1.5) rectangle (1.75,0);
			\path (1,-0.5) node[below]{$\mathcal{T}_1$};
			
			\node [empty, label = above:$u_1$] (a) at (3.5,3) {};
			\node [empty, label = above:$u_2$] (b) at (4.5,3) {};
			\node [empty, label=left:$v$] (c) at (4,2.5) {};
			\node [empty, label=left:$y$] (d) at (4,2) {};
			\node [empty, label=left:$z$] (e) at (4,1.5) {};
			\node [empty, label=above:$w$] (f) at (4.5,2) {};
			\node [empty, label=right:$s_1$] (g) at (5,2.25) {};
			\node [empty, label=right:$s_2$] (h) at (5,1.75) {};
			\node [full1] at (3.8,0.6) {};
			\node [full1] at (4.0,0.6) {};
			\node [full1] at (4.2,0.6) {};
			%\node [empty] (i) at (4,1) {};
			%\node [empty] (j) at (3.5,1) {};
			%\node [empty] (k) at (4.5,1) {};
			\draw (a) -- (c);
			\draw (b) -- (c);
			\draw (c) -- (d);
			\draw (e) -- (d);
			\draw (e) -- (f);
			\draw (g) -- (f);
			\draw (f) -- (h);
			%\draw (e) -- (i);
			%\draw (e) -- (j);
			%\draw (e) -- (k);
			\draw (e) -- (3.5,0.5);
			%\draw (i) -- (3.9,0.25);
			%\draw (i) -- (4.1,0.25);
			\draw (e) -- (4.5,0.5);
			\draw[black][dashed] (3.25,1.5) rectangle (4.75,0);
			\path (4,-0.5) node[below]{$\mathcal{T}_2$};
			
			\node [empty, label = above:$u_1$] (a) at (6.5,3) {};
			\node [empty, label = above:$u_2$] (b) at (7.5,3) {};
			\node [empty, label=left:$v$] (c) at (7,2.5) {};
			\node [empty, label=left:$y$] (d) at (7,2) {};
			\node [empty, label=left:$z$] (e) at (7,1.5) {};
			\node [empty, label = above:$s_1$] (f) at (7.5,2) {};
			\node [empty, label = above:$w$] (g) at (8,2) {};
			\node [empty, label = right:$s_2$] (h) at (8.5,2.25) {};
			\node [empty, label = right:$s_3$] (i) at (8.5,1.75) {};
			\node [full1] at (6.8,0.6) {};
			\node [full1] at (7.0,0.6) {};
			\node [full1] at (7.2,0.6) {};
			%\node [empty] (j) at (7,1) {};
			%\node [empty] (k) at (6.5,1) {};
			%\node [empty] (l) at (7.5,1) {};
			\draw (a) -- (c);
			\draw (b) -- (c);
			\draw (c) -- (d);
			\draw (e) -- (d);
			\draw (f) -- (e);
			\draw (f) -- (g);
			\draw (g) -- (h);
			\draw (i) -- (g);
			%\draw (e) -- (j);
			%\draw (e) -- (k);
			%\draw (l) -- (e);
			\draw (e) -- (6.5,0.5);
			%\draw (j) -- (6.9,0.25);
			%\draw (j) -- (7.1,0.25);
			\draw (e) -- (7.5,0.5);
			\draw[black][dashed] (6.25,1.5) rectangle (7.75,0);
			\path (7,-0.5) node[below]{$\mathcal{T}_3$};
		\end{tikzpicture}
		\caption{When $d(z)>2$, three possible families of $T$.}
		\label{case3-fig2}
	\end{figure}
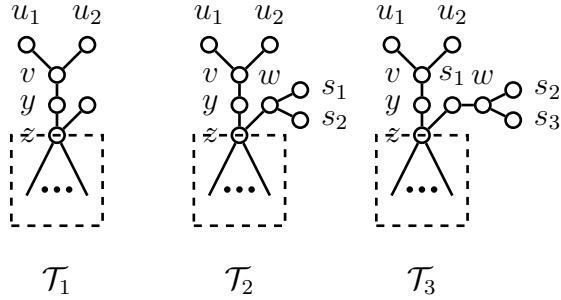

By Lemma \ref{lem1} and \ref{lem2} and the fact that $P:=u_1vyz\cdots$ is a longest path in $T$, there is at most one component in $T-z$ whose order is at least $5$,
it follows that there are only three possible families of $T$ that are pictured in Figure \ref{case3-fig2}.

$\mathcal{T}_1$: for any $T\in \mathcal{T}_1$, the vertex $z$ is adjacent to a leaf. 

$\mathcal{T}_2$: for any $T\in \mathcal{T}_2$, the vertex $z$ is adjacent to the non-leaf vertex of a $P_3$. 

$\mathcal{T}_3$: for any $T\in \mathcal{T}_3$, the vertex $z$ is adjacent to a leaf of a $K_{1,3}$. 
	
When $T\in \mathcal{T}_1$, 
\begin{align*}
	\Phi(T)=&\Phi_{\{u_1,u_2\}}(T)+2\Phi_{\{u_1,v\}}(T)+\Phi_{\{v,y\}}(T)\\
	\leq & t(n-3)+2t(n-4)+f_1(n-6)\\
	<&t(n).
\end{align*}

When $T\in\mathcal{T}_2$,
\begin{align*}
	\Phi(T)=&\Phi_{\{s_1,s_2\}}(T)+2\Phi_{\{s_1,w\}}(T)+\Phi_{\{w,z\}}(T)\\
	\leq &t(n-3)+2\Phi_{\{s_1,w\}}(T)+\Phi_{\{w,z\}}(T)\\
	\leq &t(n-3)+2\cdot4\cdot\Phi(F)+\Phi_{\{w,z\}}(T),
\end{align*}
where $F:=T-\{s_1,s_2,w,z,u_1,u_2,v,y\}$ is a forest of order $n-8$. If $F$ is a tree, then 
\[8\Phi(F)+\Phi_{\{w,z\}}(T)\leq 8t(n-8)+3f_1(n-9),\]
otherwise, 
\[8\Phi(F)+\Phi_{\{w,z\}}(T)\leq 8f_1(n-8)+3f_1(n-10).\]
Hence, 
\[\Phi(T)\leq t(n-3)+\max\{8t(n-8)+3f_1(n-9),8f_1(n-8)+3f_1(n-10)\}.\]
By direct calculations, when $n\in\{13,16,18,19\}$, 
\[\Phi(T)\leq t(n-3)+\max\{8t(n-8)+3f_1(n-9),8f_1(n-8)+3f_1(n-10)\}<t(n).\]

When $n=14$, if $F$ is not connected, then by Lemma \ref{lem1} and \ref{lem2}, $F\cong 2P_3$ and by direct calculations, $\Phi(T)<t(n)$. Now, we can assume that $F$ is a tree, we have
\[\Phi(T)\leq t(n-3)+8t(n-8)+3f_1(n-9)<t(n).\]

\begin{figure}[H]
	\begin{tikzpicture}
		[line width = 1pt, scale=0.8,
		empty/.style = {circle, draw, fill = white, inner sep=0mm, minimum size=2mm}, full/.style = {circle, draw, fill = black, inner sep=0mm, minimum size=2mm}, full1/.style = {circle, draw, fill = black, inner sep=0mm, minimum size=0.5mm}]
		
		\node [empty, label = above:$u_1$] (a) at (3.5,3) {};
		\node [empty, label = above:$u_2$] (b) at (4.5,3) {};
		\node [empty, label = left:$v$] (c) at (4,2.5) {};
		\node [empty, label = left:$y$] (d) at (4,2) {};
		\node [empty, label = left:$z$] (e) at (4,1.5) {};
		\node [empty, label = above:$w$] (f) at (4.5,2) {};
		\node [empty, label = right:$s_1$] (g) at (5,2.25) {};
		\node [empty, label = right:$s_2$] (h) at (5,1.75) {};
		\node [empty] (i) at (4,1) {};
		\node [empty] (j) at (4,0.5) {};
		\node [empty] (k) at (4,0) {};
		\node [empty] (l) at (4,-0.5) {};
		\node [empty] (m) at (4,-1) {};
		\node [empty] (n) at (3.5,-1.5) {};
		\node [empty] (o) at (4.5,-1.5) {};
		\draw (a) -- (c);
		\draw (b) -- (c);
		\draw (c) -- (d);
		\draw (e) -- (d);
		\draw (e) -- (f);
		\draw (g) -- (f);
		\draw (f) -- (h);
		\draw (e) -- (i);
		\draw (i) -- (j);
		\draw (j) -- (k);
		\draw (l) -- (k);
		\draw (l) -- (m);
		\draw (m) -- (n);
		\draw (m) -- (o);
\node [empty, label = above:$u_1$] (a) at (7.5,2.5) {};
\node [empty, label = above:$u_2$] (b) at (8.5,2.5) {};
\node [empty, label = left:$v$] (c) at (8,2) {};
\node [empty, label = left:$y$] (d) at (8,1.5) {};
\node [empty, label = left:$z$] (e) at (8,1) {};
\node [empty, label = above:$w$] (f) at (8.5,1.5) {};
\node [empty, label = right:$s_1$] (g) at (9,1.75) {};
\node [empty, label = right:$s_2$] (h) at (9,1.25) {};
\node [empty] (i) at (8.5,0.5) {};
\node [empty] (j) at (8,0.5) {};
\node [empty] (k) at (8,0) {};
\node [empty] (l) at (8,-0.5) {};
\node [empty] (m) at (8,-1) {};
\node [empty] (n) at (7.5,-1.5) {};
\node [empty] (o) at (8.5,-1.5) {};
\draw (a) -- (c);
\draw (b) -- (c);
\draw (c) -- (d);
\draw (e) -- (d);
\draw (e) -- (f);
\draw (g) -- (f);
\draw (f) -- (h);
\draw (e) -- (j);
\draw (i) -- (j);
\draw (j) -- (k);
\draw (l) -- (k);
\draw (l) -- (m);
\draw (m) -- (n);
\draw (o) -- (m);
\end{tikzpicture}
\caption{When $n=15$ and $F$ is a tree, two possible structures of $T$.}
\label{case3-fig3}
\end{figure}
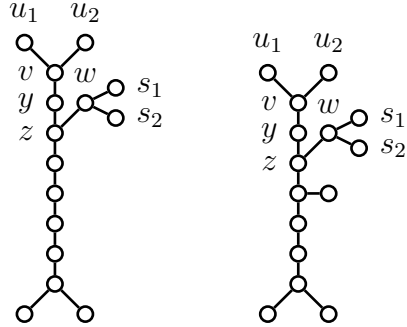

When $n=15$ and $F$ is a tree, by Lemma \ref{lem1} and \ref{lem2}, there are only two possible structures of $T$ which are pictured in Figure \ref{case3-fig3}. By direct calculations, $\Phi(T)<t(n)$. When $n=15$ and $F$ is not connected, 
\[\Phi(T)\leq t(n-3)+8f_1(n-8)+3f_1(n-10)<t(n).\]

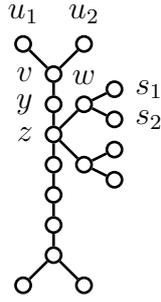
\begin{figure}[H]
	\begin{tikzpicture}
		[line width = 1pt, scale=0.8,
		empty/.style = {circle, draw, fill = white, inner sep=0mm, minimum size=2mm}, full/.style = {circle, draw, fill = black, inner sep=0mm, minimum size=2mm}, full1/.style = {circle, draw, fill = black, inner sep=0mm, minimum size=0.5mm}]
		
		\node [empty, label = above:$u_1$] (a) at (3.5,3) {};
		\node [empty, label = above:$u_2$] (b) at (4.5,3) {};
		\node [empty, label = left:$v$] (c) at (4,2.5) {};
		\node [empty, label = left:$y$] (d) at (4,2) {};
		\node [empty, label = left:$z$] (e) at (4,1.5) {};
		\node [empty, label = above:$w$] (f) at (4.5,2) {};
		\node [empty, label = right:$s_1$] (g) at (5,2.25) {};
		\node [empty, label = right:$s_2$] (h) at (5,1.75) {};
		\node [empty] (i) at (4,1) {};
		\node [empty] (j) at (4.5,1) {};
		\node [empty] (k) at (5,1.25) {};
		\node [empty] (l) at (5,0.75) {};
		\node [empty] (m) at (4,0.5) {};
		\node [empty] (n) at (4,0) {};
		\node [empty] (o) at (4,-0.5) {};
		\node [empty] (p) at (3.5,-1) {};
		\node [empty] (q) at (4.5,-1) {};
		\draw (a) -- (c);
		\draw (b) -- (c);
		\draw (c) -- (d);
		\draw (e) -- (d);
		\draw (e) -- (f);
		\draw (e) -- (j);
		\draw (f) -- (g);
		\draw (f) -- (h);
		\draw (e) -- (i);
		\draw (i) -- (m);
		\draw (m) -- (n);
		\draw (n) -- (o);
		\draw (p) -- (o);
		\draw (q) -- (o);
		\draw (j) -- (k);
		\draw (j) -- (l);
	\end{tikzpicture}
	\caption{When $n=17$ and $F$ is not connected, only one possible structure of $T$ under some assumptions.}
	\label{case3-fig4}
\end{figure}

When $n=17$ and $F$ is a tree,
\[\Phi(T)\leq t(n-3)+8t(n-8)+3f_1(n-9)<t(n).\]
When $n=17$ and $F$ is not connected, if there is a component in $F$ that is isomorphic to $K_{1,3}$, then $T$ is also belong to $\mathcal{T}_3$ and we will analysis the case in the discussion of $\mathcal{T}_3$. Now, we assume that there is no component in $F$ that is isomorphic to $K_{1,3}$. By Lemma \ref{lem1} and \ref{lem2}, there is only one possible structure of $T$ which is pictured in Figure \ref{case3-fig4}. By direct calculations, $\Phi(T)<t(n)$.

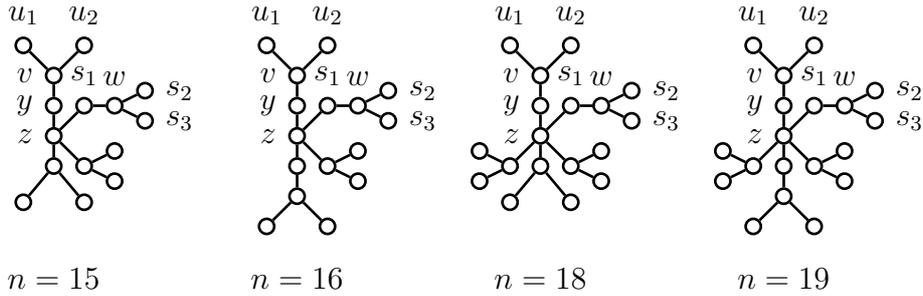
\begin{figure}[H]
	\begin{tikzpicture}
		[line width = 1pt, scale=0.8,
		empty/.style = {circle, draw, fill = white, inner sep=0mm, minimum size=2mm}, full/.style = {circle, draw, fill = black, inner sep=0mm, minimum size=2mm}, full1/.style = {circle, draw, fill = black, inner sep=0mm, minimum size=0.5mm}]
		\node [empty, label = above:$u_1$] (a) at (5.5,3) {};
		\node [empty, label = above:$u_2$] (b) at (6.5,3) {};
		\node [empty, label = left:$v$] (c) at (6,2.5) {};
		\node [empty, label = left:$y$] (d) at (6,2) {};
		\node [empty, label = left:$z$] (e) at (6,1.5) {};
		\node [empty, label = above:$s_1$] (f) at (6.5,2) {};
		\node [empty, label = above:$w$] (g) at (7,2) {};
		\node [empty, label = right:$s_2$] (h) at (7.5,2.25) {};
		\node [empty, label = right:$s_3$] (i) at (7.5,1.75) {};
		\node [empty] (j) at (6,1) {};
		\node [empty] (k) at (5.5,0.4) {};
		\node [empty] (l) at (6.5,0.4) {};
		\node [empty] (m) at (6.5,1) {};
		\node [empty] (n) at (7,1.25) {};
		\node [empty] (o) at (7,0.75) {};
		\draw (a) -- (c);
		\draw (b) -- (c);
		\draw (c) -- (d);
		\draw (e) -- (d);
		\draw (f) -- (e);
		\draw (f) -- (g);
		\draw (g) -- (h);
		\draw (i) -- (g);
		\draw (e) -- (j);
		\draw (j) -- (k);
		\draw (l) -- (j);
		\draw (m) -- (e);
		\draw (m) -- (n);
		\draw (m) -- (o);
		\path (6,-0.5) node[below]{$n=15$};
		
		\node [empty, label = above:$u_1$] (a) at (9.5,3) {};
		\node [empty, label = above:$u_2$] (b) at (10.5,3) {};
		\node [empty, label = left:$v$] (c) at (10,2.5) {};
		\node [empty, label = left:$y$] (d) at (10,2) {};
		\node [empty, label = left:$z$] (e) at (10,1.5) {};
		\node [empty, label = above:$s_1$] (f) at (10.5,2) {};
		\node [empty, label = above:$w$] (g) at (11,2) {};
		\node [empty, label = right:$s_2$] (h) at (11.5,2.25) {};
		\node [empty, label = right:$s_3$] (i) at (11.5,1.75) {};
		\node [empty] (j) at (10,1) {};
		\node [empty] (k) at (10,0.5) {};
		\node [empty] (l) at (9.5,0) {};
		\node [empty] (m) at (10.5,0) {};
		\node [empty] (n) at (10.5,1) {};
		\node [empty] (o) at (11,1.25) {};
		\node [empty] (p) at (11,0.75) {};
		\draw (a) -- (c);
		\draw (b) -- (c);
		\draw (c) -- (d);
		\draw (e) -- (d);
		\draw (f) -- (e);
		\draw (f) -- (g);
		\draw (g) -- (h);
		\draw (i) -- (g);
		\draw (e) -- (j);
		\draw (j) -- (k);
		\draw (l) -- (k);
		\draw (k) -- (m);
		\draw (o) -- (n);
		\draw (n) -- (p);
		\draw (n) -- (e);
		\path (10,-0.5) node[below]{$n=16$};
		
		\node [empty, label = above:$u_1$] (a) at (13.5,3) {};
		\node [empty, label = above:$u_2$] (b) at (14.5,3) {};
		\node [empty, label = left:$v$] (c) at (14,2.5) {};
		\node [empty, label = left:$y$] (d) at (14,2) {};
		\node [empty, label = left:$z$] (e) at (14,1.5) {};
		\node [empty, label = above:$s_1$] (f) at (14.5,2) {};
		\node [empty, label = above:$w$] (g) at (15,2) {};
		\node [empty, label = right:$s_2$] (h) at (15.5,2.25) {};
		\node [empty, label = right:$s_3$] (i) at (15.5,1.75) {};
		\node [empty] (j) at (14,1) {};
		\node [empty] (l) at (13.5,0.4) {};
		\node [empty] (m) at (14.5,0.4) {};
		\node [empty] (n) at (14.5,1) {};
		\node [empty] (o) at (15,1.25) {};
		\node [empty] (p) at (15,0.75) {};
		\node [empty] (q) at (13.5,1) {};
		\node [empty] (r) at (13,1.25) {};
		\node [empty] (s) at (13,0.75) {};
		\draw (a) -- (c);
		\draw (b) -- (c);
		\draw (c) -- (d);
		\draw (e) -- (d);
		\draw (f) -- (e);
		\draw (f) -- (g);
		\draw (g) -- (h);
		\draw (i) -- (g);
		\draw (e) -- (j);
		\draw (l) -- (j);
		\draw (j) -- (m);
		\draw (o) -- (n);
		\draw (n) -- (p);
		\draw (n) -- (e);
		\draw (q) -- (e);
		\draw (q) -- (r);
		\draw (q) -- (s);
		\path (14,-0.5) node[below]{$n=18$};
		
		\node [empty, label = above:$u_1$] (a) at (17.5,3) {};
		\node [empty, label = above:$u_2$] (b) at (18.5,3) {};
		\node [empty, label = left:$v$] (c) at (18,2.5) {};
		\node [empty, label = left:$y$] (d) at (18,2) {};
		\node [empty, label = left:$z$] (e) at (18,1.5) {};
		\node [empty, label = above:$s_1$] (f) at (18.5,2) {};
		\node [empty, label = above:$w$] (g) at (19,2) {};
		\node [empty, label = right:$s_2$] (h) at (19.5,2.25) {};
		\node [empty, label = right:$s_3$] (i) at (19.5,1.75) {};
		\node [empty] (j) at (18,1) {};
		\node [empty] (k) at (18,0.5) {};
		\node [empty] (l) at (17.5,0) {};
		\node [empty] (m) at (18.5,0) {};
		\node [empty] (n) at (18.5,1) {};
		\node [empty] (o) at (19,1.25) {};
		\node [empty] (p) at (19,0.75) {};
		\node [empty] (q) at (17.5,1) {};
		\node [empty] (r) at (17,1.25) {};
		\node [empty] (s) at (17,0.75) {};
		\draw (a) -- (c);
		\draw (b) -- (c);
		\draw (c) -- (d);
		\draw (e) -- (d);
		\draw (f) -- (e);
		\draw (f) -- (g);
		\draw (g) -- (h);
		\draw (i) -- (g);
		\draw (e) -- (j);
		\draw (j) -- (k);
		\draw (l) -- (k);
		\draw (k) -- (m);
		\draw (o) -- (n);
		\draw (n) -- (p);
		\draw (n) -- (e);
		\draw (q) -- (e);
		\draw (q) -- (r);
		\draw (q) -- (s);
		\path (18,-0.5) node[below]{$n=19$};
	\end{tikzpicture}
	\caption{When $n\in \{15,16,18,19\}$ and $\Phi(F)=f_1(n-9)$, the structures of $T$ can be completely determined.}
	\label{case3-fig5}
\end{figure}

When $T\in \mathcal{T}_3$, 
	\begin{align*}
		\Phi(T)=&\Phi_{\{u_1,u_2\}}(T)+2\Phi_{\{u_1,v\}}(T)+\Phi_{\{v,y\}}(T)\\
		\leq &t(n-3)+2t(n-4)+\Phi_{\{v,y\}}(T)\\
		=&t(n-3)+2t(n-4)+4\Phi(F),
	\end{align*}
where $F:=T-\{u_1,u_2,v,y,z,w,s_1,s_2,s_3\}$ is a forest of order $n-9$. When $n\in\{13,17\}$, 
\[\Phi(T)\leq t(n-3)+2t(n-4)+4f_1(n-9)<t(n).\]
When $n\in \{15,16,18,19\}$ and $\Phi(F)=f_1(n-9)$, by Theorem \ref{thm1} and Lemma \ref{lem1} and \ref{lem2}, we can completely determine the structures of $T$ which are pictured in Figure \ref{case3-fig5}. By direct calculations, \[\Phi(T)<t(n).\]
When $n\in \{15,16,18,19\}$ and $\Phi(F)<f_1(n-9)$, by the induction hypothesis, 
\[\Phi(T)\leq t(n-3)+2t(n-4)+4f_2(n-9)<t(n).\]
When $n=14$, by Theorem \ref{thm1} and Lemma \ref{lem1} and \ref{lem2}, it can not occur that $\Phi(F)=f_1(n-9)$, it follows that
\[\Phi(T)\leq t(n-3)+2t(n-4)+4f_2(n-9)<t(n).\]

{\bf Case 2. $d(z)=2$.}

\begin{figure}[H]
				\begin{tikzpicture}
					[line width = 1pt, scale=0.8,
					empty/.style = {circle, draw, fill = white, inner sep=0mm, minimum size=2mm}, full/.style = {circle, draw, fill = black, inner sep=0mm, minimum size=2mm}, full1/.style = {circle, draw, fill = black, inner sep=0mm, minimum size=0.5mm}]
					\node [empty] (a) at (6.75,5) {};
					\node [empty] (b) at (7.25,5) {};
					\node [empty] (c) at (7,4.5) {};
					\node [empty] (d) at (7,4) {};
					\node [empty] (e) at (7,3.5) {};
					\node [empty] (f) at (7.75,5) {};
					\node [empty] (g) at (8.25,5) {};
					\node [empty] (h) at (8,4.5) {};
					\node [empty] (i) at (8,4) {};
					\node [empty] (j) at (8,3.5) {};
					\node [empty] (k) at (7,3) {};
					\node [empty] (l) at (7.5,3) {};
					\node [empty] (m) at (8,3) {};
					\draw (a) -- (c);
					\draw (b) -- (c);
					\draw (c) -- (d);
					\draw (d) -- (e);
					\draw (f) -- (h);
					\draw (g) -- (h);
					\draw (h) -- (i);
					\draw (j) -- (i);
					\draw (k) -- (l);
					\draw (m) -- (l);
					\draw (e) -- (l);
					\draw (j) -- (l);
					
					\node [empty] (a) at (9.75,5) {};
					\node [empty] (b) at (10.25,5) {};
					\node [empty] (c) at (10,4.5) {};
					\node [empty] (d) at (10,4) {};
					\node [empty] (e) at (10,3.5) {};
					\node [empty] (f) at (10.75,5) {};
					\node [empty] (g) at (11.25,5) {};
					\node [empty] (h) at (11,4.5) {};
					\node [empty] (i) at (11,4) {};
					\node [empty] (j) at (11,3.5) {};
					\node [empty] (k) at (10,3) {};
					\node [empty] (l) at (10.5,3) {};
					\node [empty] (m) at (11,3) {};
					\draw (a) -- (c);
					\draw (b) -- (c);
					\draw (c) -- (d);
					\draw (d) -- (e);
					\draw (f) -- (h);
					\draw (g) -- (h);
					\draw (h) -- (i);
					\draw (j) -- (i);
					\draw (k) -- (l);
					\draw (m) -- (l);
					\draw (e) -- (k);
					\draw (j) -- (l);
					
					\node [empty] (a) at (12.75,5) {};
					\node [empty] (b) at (13.25,5) {};
					\node [empty] (c) at (13,4.5) {};
					\node [empty] (d) at (13,4) {};
					\node [empty] (e) at (13,3.5) {};
					\node [empty] (f) at (13.75,5) {};
					\node [empty] (g) at (14.25,5) {};
					\node [empty] (h) at (14,4.5) {};
					\node [empty] (i) at (14,4) {};
					\node [empty] (j) at (14,3.5) {};
					\node [empty] (k) at (13,3) {};
					\node [empty] (l) at (13.5,3) {};
					\node [empty] (m) at (14,3) {};
					\draw (a) -- (c);
					\draw (b) -- (c);
					\draw (c) -- (d);
					\draw (d) -- (e);
					\draw (f) -- (h);
					\draw (g) -- (h);
					\draw (h) -- (i);
					\draw (j) -- (i);
					\draw (k) -- (l);
					\draw (m) -- (l);
					\draw (e) -- (k);
					\draw (j) -- (m);
				\end{tikzpicture}
				\caption{When $n=13$ and $d(z)=2$, three possible structures of $T$.}
				\label{case3-fig6}
			\end{figure}
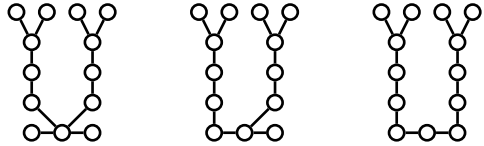

When $n=13$, by Lemma \ref{lem1} and \ref{lem2}, there are only three possible structures of $T$ which are pictured in Figure \ref{case3-fig6}. By direct calculations, 
\[\Phi(T)<t(n).\]

\begin{figure}[H]
				\begin{tikzpicture}
					[line width = 1pt, scale=0.8,
					empty/.style = {circle, draw, fill = white, inner sep=0mm, minimum size=2mm}, full/.style = {circle, draw, fill = black, inner sep=0mm, minimum size=2mm}, full1/.style = {circle, draw, fill = black, inner sep=0mm, minimum size=0.5mm}]
					\node [empty] (a) at (6.75,5) {};
					\node [empty] (b) at (7.25,5) {};
					\node [empty] (c) at (7,4.5) {};
					\node [empty] (d) at (7,4) {};
					\node [empty] (e) at (7,3.5) {};
					\node [empty] (f) at (7.75,5) {};
					\node [empty] (g) at (8.25,5) {};
					\node [empty] (h) at (8,4.5) {};
					\node [empty] (i) at (8,4) {};
					\node [empty] (j) at (8,3.5) {};
					\node [empty] (k) at (7,3) {};
					\node [empty] (l) at (7.5,3) {};
					\node [empty] (m) at (8,3) {};
					\node [empty] (n) at (7.5,2.5) {};
					\draw (a) -- (c);
					\draw (b) -- (c);
					\draw (c) -- (d);
					\draw (d) -- (e);
					\draw (f) -- (h);
					\draw (g) -- (h);
					\draw (h) -- (i);
					\draw (j) -- (i);
					\draw (k) -- (l);
					\draw (m) -- (l);
					\draw (e) -- (k);
					\draw (j) -- (k);
					\draw (n) -- (l);
					
					\node [empty] (a) at (9.75,5) {};
					\node [empty] (b) at (10.25,5) {};
					\node [empty] (c) at (10,4.5) {};
					\node [empty] (d) at (10,4) {};
					\node [empty] (e) at (10,3.5) {};
					\node [empty] (f) at (10.75,5) {};
					\node [empty] (g) at (11.25,5) {};
					\node [empty] (h) at (11,4.5) {};
					\node [empty] (i) at (11,4) {};
					\node [empty] (j) at (11,3.5) {};
					\node [empty] (k) at (10,3) {};
					\node [empty] (l) at (10.5,3) {};
					\node [empty] (m) at (11,3) {};
					\node [empty] (n) at (10.5,2.5) {};
					\draw (a) -- (c);
					\draw (b) -- (c);
					\draw (c) -- (d);
					\draw (d) -- (e);
					\draw (f) -- (h);
					\draw (g) -- (h);
					\draw (h) -- (i);
					\draw (j) -- (i);
					\draw (k) -- (l);
					\draw (m) -- (l);
					\draw (e) -- (k);
					\draw (j) -- (l);
					\draw (n) -- (l);
					
					\node [empty] (a) at (12.75,5) {};
					\node [empty] (b) at (13.25,5) {};
					\node [empty] (c) at (13,4.5) {};
					\node [empty] (d) at (13,4) {};
					\node [empty] (e) at (13,3.5) {};
					\node [empty] (f) at (13.75,5) {};
					\node [empty] (g) at (14.25,5) {};
					\node [empty] (h) at (14,4.5) {};
					\node [empty] (i) at (14,4) {};
					\node [empty] (j) at (14,3.5) {};
					\node [empty] (k) at (13,3) {};
					\node [empty] (l) at (13.5,3) {};
					\node [empty] (m) at (14,3) {};
					\node [empty] (n) at (13.5,2.5) {};
					\draw (a) -- (c);
					\draw (b) -- (c);
					\draw (c) -- (d);
					\draw (d) -- (e);
					\draw (f) -- (h);
					\draw (g) -- (h);
					\draw (h) -- (i);
					\draw (j) -- (i);
					\draw (k) -- (l);
					\draw (m) -- (l);
					\draw (e) -- (k);
					\draw (j) -- (m);
					\draw (n) -- (l);
					
					\node [empty] (a) at (6.5,1) {};
					\node [empty] (b) at (7,1) {};
					\node [empty] (c) at (6.75,0.5) {};
					\node [empty] (d) at (6.75,0) {};
					\node [empty] (e) at (6.75,-0.5) {};
					\node [empty] (f) at (8,1) {};
					\node [empty] (g) at (8.5,1) {};
					\node [empty] (h) at (8.25,0.5) {};
					\node [empty] (i) at (8.25,0) {};
					\node [empty] (j) at (8.25,-0.5) {};
					\node [empty] (k) at (7.25,-1) {};
					\node [empty] (l) at (6.75,-1) {};
					\node [empty] (m) at (7.75,-1) {};
					\node [empty] (n) at (8.25,-1) {};
					\draw (a) -- (c);
					\draw (b) -- (c);
					\draw (c) -- (d);
					\draw (d) -- (e);
					\draw (f) -- (h);
					\draw (g) -- (h);
					\draw (h) -- (i);
					\draw (j) -- (i);
					\draw (k) -- (l);
					\draw (m) -- (k);
					\draw (m) -- (n);
					\draw (e) -- (l);
					\draw (j) -- (m);
					
					\node [empty] (a) at (9.5,1) {};
					\node [empty] (b) at (10,1) {};
					\node [empty] (c) at (9.75,0.5) {};
					\node [empty] (d) at (9.75,0) {};
					\node [empty] (e) at (9.75,-0.5) {};
					\node [empty] (f) at (11,1) {};
					\node [empty] (g) at (11.5,1) {};
					\node [empty] (h) at (11.25,0.5) {};
					\node [empty] (i) at (11.25,0) {};
					\node [empty] (j) at (11.25,-0.5) {};
					\node [empty] (k) at (10.25,-1) {};
					\node [empty] (l) at (9.75,-1) {};
					\node [empty] (m) at (10.75,-1) {};
					\node [empty] (n) at (11.25,-1) {};
					\draw (a) -- (c);
					\draw (b) -- (c);
					\draw (c) -- (d);
					\draw (d) -- (e);
					\draw (f) -- (h);
					\draw (g) -- (h);
					\draw (h) -- (i);
					\draw (j) -- (i);
					\draw (k) -- (l);
					\draw (m) -- (k);
					\draw (m) -- (n);
					\draw (e) -- (k);
					\draw (j) -- (m);
					
					\node [empty] (a) at (12.5,1) {};
					\node [empty] (b) at (13,1) {};
					\node [empty] (c) at (12.75,0.5) {};
					\node [empty] (d) at (12.75,0) {};
					\node [empty] (e) at (12.75,-0.5) {};
					\node [empty] (f) at (14,1) {};
					\node [empty] (g) at (14.5,1) {};
					\node [empty] (h) at (14.25,0.5) {};
					\node [empty] (i) at (14.25,0) {};
					\node [empty] (j) at (14.25,-0.5) {};
					\node [empty] (k) at (13.25,-1) {};
					\node [empty] (l) at (12.75,-1) {};
					\node [empty] (m) at (13.75,-1) {};
					\node [empty] (n) at (14.25,-1) {};
					\draw (a) -- (c);
					\draw (b) -- (c);
					\draw (c) -- (d);
					\draw (d) -- (e);
					\draw (f) -- (h);
					\draw (g) -- (h);
					\draw (h) -- (i);
					\draw (j) -- (i);
					\draw (k) -- (l);
					\draw (m) -- (k);
					\draw (m) -- (n);
					\draw (e) -- (l);
					\draw (j) -- (n);
				\end{tikzpicture}
				\caption{When $n=14$ and $d(z)=2$, six possible structures of $T$.}
				\label{case3-fig7}
			\end{figure}
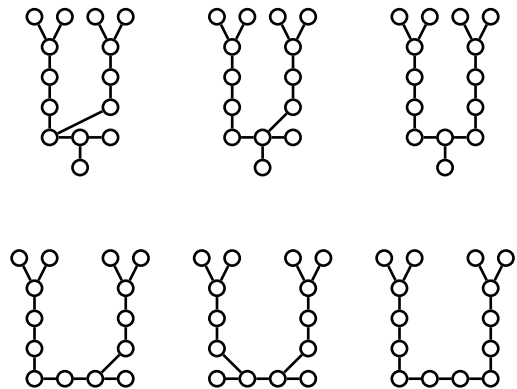

When $n=14$, by Lemma \ref{lem1} and \ref{lem2}, there are only six possible structures of $T$ which are pictured in Figure \ref{case3-fig7}. By direct calculations, 
\[\Phi(T)<t(n).\]

\begin{figure}[H]
				\begin{tikzpicture}
					[line width = 1pt, scale=0.8,
					empty/.style = {circle, draw, fill = white, inner sep=0mm, minimum size=2mm}, full/.style = {circle, draw, fill = black, inner sep=0mm, minimum size=2mm}, full1/.style = {circle, draw, fill = black, inner sep=0mm, minimum size=0.5mm}]
					
					\node [empty] (a) at (6.5,1) {};
					\node [empty] (b) at (7,1) {};
					\node [empty] (c) at (6.75,0.5) {};
					\node [empty] (d) at (6.75,0) {};
					\node [empty] (e) at (6.75,-0.5) {};
					\node [empty] (f) at (8,1) {};
					\node [empty] (g) at (8.5,1) {};
					\node [empty] (h) at (8.25,0.5) {};
					\node [empty] (i) at (8.25,0) {};
					\node [empty] (j) at (8.25,-0.5) {};
					\node [empty] (k) at (7.25,-1) {};
					\node [empty] (l) at (6.75,-1) {};
					\node [empty] (m) at (7.75,-1) {};
					\node [empty] (n) at (8.25,-1) {};
					\node [empty] (o) at (7.25,-1.5) {};
					\draw (a) -- (c);
					\draw (b) -- (c);
					\draw (c) -- (d);
					\draw (d) -- (e);
					\draw (f) -- (h);
					\draw (g) -- (h);
					\draw (h) -- (i);
					\draw (j) -- (i);
					\draw (k) -- (l);
					\draw (m) -- (k);
					\draw (m) -- (n);
					\draw (e) -- (n);
					\draw (j) -- (n);
					\draw (o) -- (k);
					
					\node [empty] (a) at (9.5,1) {};
					\node [empty] (b) at (10,1) {};
					\node [empty] (c) at (9.75,0.5) {};
					\node [empty] (d) at (9.75,0) {};
					\node [empty] (e) at (9.75,-0.5) {};
					\node [empty] (f) at (11,1) {};
					\node [empty] (g) at (11.5,1) {};
					\node [empty] (h) at (11.25,0.5) {};
					\node [empty] (i) at (11.25,0) {};
					\node [empty] (j) at (11.25,-0.5) {};
					\node [empty] (k) at (10.25,-1) {};
					\node [empty] (l) at (9.75,-1) {};
					\node [empty] (m) at (10.75,-1) {};
					\node [empty] (n) at (11.25,-1) {};
					\node [empty] (o) at (10.25,-1.5) {};
					\draw (a) -- (c);
					\draw (b) -- (c);
					\draw (c) -- (d);
					\draw (d) -- (e);
					\draw (f) -- (h);
					\draw (g) -- (h);
					\draw (h) -- (i);
					\draw (j) -- (i);
					\draw (k) -- (l);
					\draw (m) -- (k);
					\draw (m) -- (n);
					\draw (e) -- (k);
					\draw (j) -- (m);
					\draw (o) -- (k);
					
					\node [empty] (a) at (12.5,1) {};
					\node [empty] (b) at (13,1) {};
					\node [empty] (c) at (12.75,0.5) {};
					\node [empty] (d) at (12.75,0) {};
					\node [empty] (e) at (12.75,-0.5) {};
					\node [empty] (f) at (14,1) {};
					\node [empty] (g) at (14.5,1) {};
					\node [empty] (h) at (14.25,0.5) {};
					\node [empty] (i) at (14.25,0) {};
					\node [empty] (j) at (14.25,-0.5) {};
					\node [empty] (k) at (13.25,-1) {};
					\node [empty] (l) at (12.75,-1) {};
					\node [empty] (m) at (13.75,-1) {};
					\node [empty] (n) at (14.25,-1) {};
					\node [empty] (o) at (13.25,-1.5) {};
					\draw (a) -- (c);
					\draw (b) -- (c);
					\draw (c) -- (d);
					\draw (d) -- (e);
					\draw (f) -- (h);
					\draw (g) -- (h);
					\draw (h) -- (i);
					\draw (j) -- (i);
					\draw (k) -- (l);
					\draw (m) -- (k);
					\draw (m) -- (n);
					\draw (e) -- (m);
					\draw (j) -- (n);
					\draw (o) -- (k);
					
					\node [empty] (a) at (15.5,1) {};
					\node [empty] (b) at (16,1) {};
					\node [empty] (c) at (15.75,0.5) {};
					\node [empty] (d) at (15.75,0) {};
					\node [empty] (e) at (15.75,-0.5) {};
					\node [empty] (f) at (17,1) {};
					\node [empty] (g) at (17.5,1) {};
					\node [empty] (h) at (17.25,0.5) {};
					\node [empty] (i) at (17.25,0) {};
					\node [empty] (j) at (17.25,-0.5) {};
					\node [empty] (k) at (15.75,-1) {};
					\node [empty] (l) at (16.125,-1) {};
					\node [empty] (m) at (16.5,-1) {};
					\node [empty] (n) at (16.875,-1) {};
					\node [empty] (o) at (17.25,-1) {};
					\draw (a) -- (c);
					\draw (b) -- (c);
					\draw (c) -- (d);
					\draw (d) -- (e);
					\draw (f) -- (h);
					\draw (g) -- (h);
					\draw (h) -- (i);
					\draw (j) -- (i);
					\draw (k) -- (l);
					\draw (l) -- (m);
					\draw (n) -- (m);
					\draw (n) -- (o);
					\draw (e) -- (l);
					\draw (j) -- (n);
					
					\node [empty] (a) at (18.5,1) {};
					\node [empty] (b) at (19,1) {};
					\node [empty] (c) at (18.75,0.5) {};
					\node [empty] (d) at (18.75,0) {};
					\node [empty] (e) at (18.75,-0.5) {};
					\node [empty] (f) at (20,1) {};
					\node [empty] (g) at (20.5,1) {};
					\node [empty] (h) at (20.25,0.5) {};
					\node [empty] (i) at (20.25,0) {};
					\node [empty] (j) at (20.25,-0.5) {};
					\node [empty] (k) at (18.75,-1) {};
					\node [empty] (l) at (19.125,-1) {};
					\node [empty] (m) at (19.5,-1) {};
					\node [empty] (n) at (19.875,-1) {};
					\node [empty] (o) at (20.25,-1) {};
					\draw (a) -- (c);
					\draw (b) -- (c);
					\draw (c) -- (d);
					\draw (d) -- (e);
					\draw (f) -- (h);
					\draw (g) -- (h);
					\draw (h) -- (i);
					\draw (j) -- (i);
					\draw (k) -- (l);
					\draw (l) -- (m);
					\draw (n) -- (m);
					\draw (n) -- (o);
					\draw (e) -- (k);
					\draw (j) -- (n);
					
					\node [empty] (a) at (6.5,-3) {};
					\node [empty] (b) at (7,-3) {};
					\node [empty] (c) at (6.75,-3.5) {};
					\node [empty] (d) at (6.75,-4) {};
					\node [empty] (e) at (6.75,-4.5) {};
					\node [empty] (f) at (8,-3) {};
					\node [empty] (g) at (8.5,-3) {};
					\node [empty] (h) at (8.25,-3.5) {};
					\node [empty] (i) at (8.25,-4) {};
					\node [empty] (j) at (8.25,-4.5) {};
					\node [empty] (k) at (7.25,-5) {};
					\node [empty] (l) at (6.75,-5) {};
					\node [empty] (m) at (7.75,-5) {};
					\node [empty] (n) at (8.25,-5) {};
					\node [empty] (o) at (7.25,-5.5) {};
					\draw (a) -- (c);
					\draw (b) -- (c);
					\draw (c) -- (d);
					\draw (d) -- (e);
					\draw (f) -- (h);
					\draw (g) -- (h);
					\draw (h) -- (i);
					\draw (j) -- (i);
					\draw (k) -- (l);
					\draw (m) -- (k);
					\draw (m) -- (n);
					\draw (e) -- (l);
					\draw (j) -- (m);
					\draw (o) -- (k);
					
					\node [empty] (a) at (9.5,-3) {};
					\node [empty] (b) at (10,-3) {};
					\node [empty] (c) at (9.75,-3.5) {};
					\node [empty] (d) at (9.75,-4) {};
					\node [empty] (e) at (9.75,-4.5) {};
					\node [empty] (f) at (11,-3) {};
					\node [empty] (g) at (11.5,-3) {};
					\node [empty] (h) at (11.25,-3.5) {};
					\node [empty] (i) at (11.25,-4) {};
					\node [empty] (j) at (11.25,-4.5) {};
					\node [empty] (k) at (10.25,-5) {};
					\node [empty] (l) at (9.75,-5) {};
					\node [empty] (m) at (10.75,-5) {};
					\node [empty] (n) at (11.25,-5) {};
					\node [empty] (o) at (10.25,-5.5) {};
					\draw (a) -- (c);
					\draw (b) -- (c);
					\draw (c) -- (d);
					\draw (d) -- (e);
					\draw (f) -- (h);
					\draw (g) -- (h);
					\draw (h) -- (i);
					\draw (j) -- (i);
					\draw (k) -- (l);
					\draw (m) -- (k);
					\draw (m) -- (n);
					\draw (e) -- (k);
					\draw (j) -- (n);
					\draw (o) -- (k);
					
					\node [empty] (a) at (12.5,-3) {};
					\node [empty] (b) at (13,-3) {};
					\node [empty] (c) at (12.75,-3.5) {};
					\node [empty] (d) at (12.75,-4) {};
					\node [empty] (e) at (12.75,-4.5) {};
					\node [empty] (f) at (14,-3) {};
					\node [empty] (g) at (14.5,-3) {};
					\node [empty] (h) at (14.25,-3.5) {};
					\node [empty] (i) at (14.25,-4) {};
					\node [empty] (j) at (14.25,-4.5) {};
					\node [empty] (k) at (13.25,-5) {};
					\node [empty] (l) at (12.75,-5) {};
					\node [empty] (m) at (13.75,-5) {};
					\node [empty] (n) at (14.25,-5) {};
					\node [empty] (o) at (13.25,-5.5) {};
					\draw (a) -- (c);
					\draw (b) -- (c);
					\draw (c) -- (d);
					\draw (d) -- (e);
					\draw (f) -- (h);
					\draw (g) -- (h);
					\draw (h) -- (i);
					\draw (j) -- (i);
					\draw (k) -- (l);
					\draw (m) -- (k);
					\draw (m) -- (n);
					\draw (e) -- (l);
					\draw (j) -- (n);
					\draw (o) -- (k);
					
					\node [empty] (a) at (15.5,-3) {};
					\node [empty] (b) at (16,-3) {};
					\node [empty] (c) at (15.75,-3.5) {};
					\node [empty] (d) at (15.75,-4) {};
					\node [empty] (e) at (15.75,-4.5) {};
					\node [empty] (f) at (17,-3) {};
					\node [empty] (g) at (17.5,-3) {};
					\node [empty] (h) at (17.25,-3.5) {};
					\node [empty] (i) at (17.25,-4) {};
					\node [empty] (j) at (17.25,-4.5) {};
					\node [empty] (k) at (16.5,-5) {};
					\node [empty] (l) at (15.75,-5) {};
					\node [empty] (m) at (16.75,-5.5) {};
					\node [empty] (n) at (17.25,-5) {};
					\node [empty] (o) at (16.25,-5.5) {};
					\draw (a) -- (c);
					\draw (b) -- (c);
					\draw (c) -- (d);
					\draw (d) -- (e);
					\draw (f) -- (h);
					\draw (g) -- (h);
					\draw (h) -- (i);
					\draw (j) -- (i);
					\draw (k) -- (l);
					\draw (m) -- (k);
					\draw (k) -- (n);
					\draw (e) -- (l);
					\draw (j) -- (n);
					\draw (o) -- (k);
					
					\node [empty] (a) at (18.5,-3) {};
					\node [empty] (b) at (19,-3) {};
					\node [empty] (c) at (18.75,-3.5) {};
					\node [empty] (d) at (18.75,-4) {};
					\node [empty] (e) at (18.75,-4.5) {};
					\node [empty] (f) at (20,-3) {};
					\node [empty] (g) at (20.5,-3) {};
					\node [empty] (h) at (20.25,-3.5) {};
					\node [empty] (i) at (20.25,-4) {};
					\node [empty] (j) at (20.25,-4.5) {};
					\node [empty] (k) at (18.75,-5) {};
					\node [empty] (l) at (19.125,-5) {};
					\node [empty] (m) at (19.5,-5) {};
					\node [empty] (n) at (19.875,-5) {};
					\node [empty] (o) at (20.25,-5) {};
					\draw (a) -- (c);
					\draw (b) -- (c);
					\draw (c) -- (d);
					\draw (d) -- (e);
					\draw (f) -- (h);
					\draw (g) -- (h);
					\draw (h) -- (i);
					\draw (j) -- (i);
					\draw (k) -- (l);
					\draw (l) -- (m);
					\draw (n) -- (m);
					\draw (n) -- (o);
					\draw (e) -- (k);
					\draw (j) -- (o);
				\end{tikzpicture}
				\caption{When $n=15$ and $d(z)=2$, ten possible structures of $T$.}
				\label{case3-fig8}
			\end{figure}
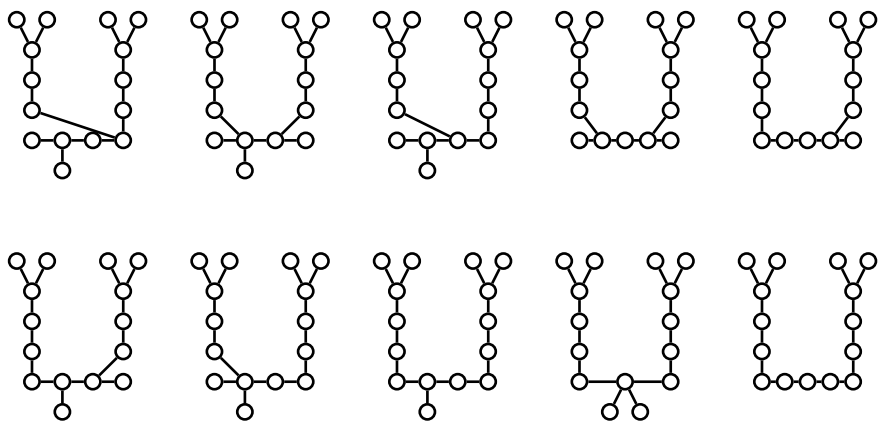

When $n=15$, by Lemma \ref{lem1} and \ref{lem2}, there are only ten possible structures of $T$ which are pictured in Figure \ref{case3-fig8}. By direct calculations, 
\[\Phi(T)<t(n).\]

When $16\leq n\leq 19,$
\begin{align*}
	\Phi(T)=&\Phi_{\{u_1,u_2\}}(T)+2\Phi_{\{u_1,v\}}(T)+\Phi_{\{v,y\}}(T)\\
	\leq &t(n-3)+2t(n-4)+t(n-5)\\
	<&t(n).
\end{align*}

In conclusion, we have proved that when $11\leq n\leq 19$, for any tree $T$ of order $n$, 
	\[\Phi(T)\leq t(n),\]
with equality if and only if $n\equiv 1\pmod3$ and $T\cong T^*_n$. \qed

\section*{Declaration of competing interest}

The authors have no relevant financial or non-financial interests to disclose.

\section*{Data availability}

No data was used for the research described in the article.

\section*{Acknowledgements}

This work was supported by Beijing Natural Science Foundation (No. 1232005).

	\bibliographystyle{unsrt}
	
\end{document}